\documentclass[11pt]{article}
\usepackage[letterpaper,left=1in,right=1in,top=1in,bottom=1in]{geometry}
\usepackage[mathscr]{eucal}
\usepackage{epsfig,epsf,psfrag}
\usepackage{amssymb,amsfonts,latexsym}
\usepackage{amsmath}
\usepackage{graphicx}
\usepackage{epstopdf}
\usepackage{bm,xcolor,url}
\usepackage{fixltx2e}%ordering of single and double column floats
\usepackage{array}%array and tabular environments
\usepackage{verbatim}
\usepackage[noend]{algpseudocode}
\usepackage{cite}
\usepackage{algorithm}
\usepackage{verbatim}
\usepackage{textcomp}
\usepackage{mathrsfs}
\usepackage[referable]{threeparttablex}
\usepackage{subfig}
\usepackage{amsthm}
\usepackage{enumerate}
\usepackage{footnote}
\usepackage{enumitem}
\usepackage{xr}
\usepackage{mathtools}
\catcode`~=11 \def\UrlSpecials{\do\~{\kern -.15em\lower .7ex\hbox{~}\kern .04em}} \catcode`~=13 

\allowdisplaybreaks[3]

\newcommand{\tnorm}[1]{{\left\vert\kern-0.25ex\left\vert\kern-0.25ex\left\vert #1 
    \right\vert\kern-0.25ex\right\vert\kern-0.25ex\right\vert}}
\newcommand{\tnormt}[1]{{\vert\kern-0.25ex\vert\kern-0.25ex\vert #1 
    \vert\kern-0.25ex\vert\kern-0.25ex\vert}}

\newcommand{\norm}[1]{\left\Vert#1\right\Vert}
\newcommand{\normt}[1]{\Vert#1\Vert}

\newcommand{\nn}{\nonumber}

\newcommand{\prox}{\mathsf{prox}}
\newcommand{\dom}{\mathsf{dom}\,}

\newcommand{\where}{\mathrm{where}}

\newcommand{\andd}{\mathrm{and}}

\newcommand{\conv}{\mathsf{conv}\,}

\newcommand{\calA}{\mathcal{A}}

\newcommand{\calV}{\mathcal{V}}

\newcommand{\calX}{\mathcal{X}}

\newcommand{\rmP}{\mathrm{P}}

\newcommand{\bbR}{\mathbb{R}}

\newcommand{\barbbR}{\bar{\bbR}}

\DeclareMathAlphabet{\mathbsf}{OT1}{cmss}{bx}{n}
\newcommand{\barg}{\bar{g}}

\newcommand{\barx}{\bar{x}}
\newcommand{\bary}{\bar{y}}

\newcommand{\barQ}{\bar{Q}}

\newcommand{\uG}{\underline{G}}

\newcommand{\uJ}{\underline{J}}

\newcommand{\ip}[2]{\left\langle{#1},{#2}\right\rangle}
\newcommand{\ipt}[2]{\langle{#1},{#2}\rangle}

\DeclareMathOperator*{\argmin}{arg\,min}

\newtheorem{theorem}{Theorem} 
\newtheorem*{theorem*}{Theorem}
\newtheorem{lemma}{Lemma}

\newtheorem{prop}{Proposition}
\newtheorem{corollary}{Corollary}
 
\newtheorem*{assump*}{Assumption}

\theoremstyle{definition}
\newtheorem{definition}{Definition}

\theoremstyle{remark}

\newcommand{\qednew}{\nobreak \ifvmode \relax \else
      \ifdim\lastskip<1.5em \hskip-\lastskip
      \hskip1.5em plus0em minus0.5em \fi \nobreak
      \vrule height0.75em width0.5em depth0.25em\fi}

\title{An Optimal First-Order Active-Hull Method for Convex Composite 
 H\"older Smooth Optimization} % Part I: The Bounded-Domain Case} % With First-Order and Generalized 
\usepackage[colorlinks=true,breaklinks=true,bookmarks=true,urlcolor=blue,
     citecolor=blue,linkcolor=blue,bookmarksopen=false,draft=false]{hyperref}
\author{Renbo Zhao
\thanks{Tippie College of Business, University of Iowa, Iowa City, IA 52242 ({mailto:  renbo-zhao@uiowa.edu}).}
}
\usepackage[normalem]{ulem}

\usepackage{algorithm}
\newcommand{\gap}{{\mathsf{gap}}}
\newcommand{\inp}{{\rm in}}
\newcommand{\FW}{\mathrm{FW}}
\newcommand{\TA}{\mathrm{TA}}

\numberwithin{equation}{section}
\numberwithin{definition}{section}
\numberwithin{prop}{section}
\numberwithin{lemma}{section}
\numberwithin{theorem}{section}
\numberwithin{remark}{section}
\numberwithin{corollary}{section}

\begin{document}

\maketitle

\begin{abstract}
We consider the problem  $(\rmP): {\min}_{x\in \mathbb{R}^n}\, f(x)+g(x)$ under the Euclidean geometry, where $f$ is convex and $(M,\nu)$-H\"older smooth on $\mathbb{R}^n$, and $g$ is proper, closed and convex on $\mathbb{R}^n$,  admits a generalized linear minimization oracle (GLMO), and has bounded  domain with diameter $D$. An first-order active-hull  (AH) method has access to both the first-order oracle (FO) of $f$ and the GLMO of $g$, and outputs a point that lies in the convex hull of all the past GLMO returns. % Under the Euclidean geometry (and in the high-dimensional regime),  
We   propose a new first-order AH method for solving~(P), % to obtain an $\epsilon$-optimal solution of~(P), 
which is based on our newly developed ``recentered homotopy smoothing framework''. We show that in the high-dimensional regime, this method achieves the optimal FO and GLMO complexties %in~\eqref{eq:opt_AH_nu} 
{simultaneously for all} $\nu\in[0,1]$. %, which are 
%\begin{align}
%O\Big( M^{\frac{2}{1+3\nu}} D^{\frac{2(1+\nu)}{1+3\nu}} \epsilon^{-\frac{2}{1+3\nu}}
%\Big)\quad\andd\quad  O\Big(M^{\frac{2}{1+\nu}}D^2\epsilon^{-\frac{2}{1+\nu}}\Big), \quad \mbox{respectively}, \quad \mbox{for all }\, \nu\in[0,1]. \nn
%%\label{eq:opt_AH_nu}
%\end{align}
 %$O\big( M^{\frac{2}{1+3\nu}} D^{\frac{2(1+\nu)}{1+3\nu}} \epsilon^{-\frac{2}{1+3\nu}}\big)$ and  $O\big(M^{\frac{2}{1+\nu}}D^2\epsilon^{-\frac{2}{1+\nu}}\big)$, \mbox{respectively},
\end{abstract}

\section{Introduction} \label{sec:intro}

Consider the following convex composite optimization problem: 
\begin{equation}
P^*:={\min}_{x\in \bbR^n}\; \{P(y):=f(x)+g(x)\}, \tag{P} \label{eq:P}
\end{equation}
where %$A:\bbR^n\to \bbR^m$ is a linear operator and 
$f:\bbR^n\to \bbR$ is convex and $(M,\nu)$-H\"older smooth on $\bbR^n$ for some $M>0$ and $\nu\in[0,1]$, i.e., 
\begin{equation}
0\le f(x)- f(z) - \ip{g_z}{x-z} \le \frac{M}{1+\nu}\norm{x-z}^{1+\nu}, \quad \forall\,g_z\in \partial f(z),\quad \forall\, z,x\in \bbR^m. \label{eq:holder}
\end{equation}
%where 
Throughout this paper, $\normt{\cdot}$ %denotes some norm on $\bbR^m$. %
denotes  the {\em Euclidean} norm. In particular, we call $f$ %{\em nonsmooth} (and 
$M$-Lipschitz if $\nu=0$ and {$M$-smooth} if $\nu=1$. 
Also, $g:\bbR^n\to \barbbR:=\bbR\cup\{+\infty\}$ is a proper, closed and convex (p.c.c.) function.  Additionally, we make the following two blanket assumptions: 
\begin{enumerate}[label=(A\arabic*)]
%\item \label{assump:f_holder} At any $x\in\R^m$, $f(x)$ and a subgradient of $f$, i.e., $g_x\in\partial f(x)$, can be easily computed.  The oracle that outputs $(f(x),g_x)$ at any $x\in\R^m$ 
%is called the first-order oracle of $f$. 
\item \label{assump:LMO} For any $\lambda\in \bbR^n$, the following optimization problem % and there is an oracle for the mapping 
\begin{equation}
 {\min}_{x\in \bbR^n}\,\{\ip{\lambda}{y}+g(y)\}, %\tag{LM}
  \label{eq:LMO}
\end{equation}
admits an easily computable optimal solution. We shall call the oracle that outputs an solution of~\eqref{eq:LMO} the   {\em generalized linear minimization oracle} (GLMO) of $g$.
\item \label{assump:domg_bdd} $\dom g$ is bounded, where $\dom g:= \{y\in\bbR^n:g(y)<+\infty\}$ denotes the domain of $g$. Consequently, define
\begin{equation}
D:= {\sup}_{x,\,z\in \dom g}\;\normt{x-z}<+\infty.  \label{eq:def_D}
\end{equation}
\end{enumerate}
Denote the set of optimal solutions of~\eqref{eq:P} by $\calX^*$. Clearly, under~\ref{assump:domg_bdd}, we know that $\calX^*\ne \emptyset.$
% and  

In this work, we aim to investigate the computational complexities of the class of active-hull methods for solving~\eqref{eq:P} with access to the first-order oracle (FO) of $f$ and the GLMO of $g$.  %first-order % investigates the 
The class of active-hull (AH) methods was initially proposed in~\cite{Lan_13} when $g=\iota_\calX$ (i.e., %namely 
the indicator function of some convex and compact set $\calX$), % that is convex and compact), 
and generalized to accommodate any p.c.c.\ function $g$ in~\cite{Zhao_26a}. This class of methods, as shown in Algorithm~\ref{algo:LMO-AH}, has a extremely simple structure.  % and is shown in 
In essence, the output $\bary_K$ has to lie in the convex hull of all the past output of the GLMO %, together with 
and the starting point $y_0$. As mentioned in~\cite{Lan_13,Zhao_26a}, this class of methods is fairly broad, and includes the Frank-Wolfe (FW) method~\cite{Frank_56} and many of its variants (see~\cite{Lan_13,Jaggi_13,Lacoste_15,Lan_16,Zhao_22} and references therein) as special cases. In~\cite{Zhao_26a}, it was shown that in the high-dimensional regime, the optimal GLMO complexity of % Algorithm~\ref{algo:LMO-AH} 
the class of AH methods to obtain an $\epsilon$-optimal solution of~\eqref{eq:P} %. (Here we call 
(i.e., a point $y\in \dom g$ such that %$\epsilon$-optimal solution of~\eqref{eq:P} if 
$P(y) - P^*<\epsilon$)  is %for solving~\eqref{eq:P} is 
\begin{align}
O\left(M^{\frac{2}{1+\nu}}D^2\epsilon^{-\frac{2}{1+\nu}}\right), \quad \forall\,\nu\in[0,1].  \label{eq:opt_GLMO_AH_nu}
 \end{align}

Although the class of AH methods does not specify the oracle access of $f$, we are particularly interested in the subclass of methods with FO access of $f$, which we call {\em the first-order AH methods}. Since the first-order AH methods use both the FO of $f$ and GLMO of $g$ to solve~\eqref{eq:P}, 
a natural question to ask is %investigate is 
%one may naturally wonder what 
the optimal complexity of either oracle %should be, in order 
for a first-order AH method to obtain an $\epsilon$-optimal solution of~\eqref{eq:P}. %(Here we call $y\in \dom g$ an $\epsilon$-optimal solution of~\eqref{eq:P} if $P(y) - P^*<\epsilon$.)  
In fact, this question has been resolved together in~\cite{Ouyang_23} and~\cite{Zhao_26a}, which proposed and analyzed two different first-order AH methods, respectively. For convenience, let $N_{\rm FO}(\epsilon)$ and $N_{\rm GLMO}(\epsilon)$ denote the numbers of FO and GLMO calls made by any first-order AH method
%, respectively,  
to obtain  an $\epsilon$-optimal solution of~\eqref{eq:P}. In~\cite{Ouyang_23}, the universal conditional gradient sliding (CGS) method was proposed, and the authors showed that for all $ \nu\in[0,1],$ this method achieves the following FO and GLMO complexities: 
\begin{align}
N_{\rm FO} (\epsilon)= O\Big( M^{\frac{2}{1+3\nu}} D^{\frac{2(1+\nu)}{1+3\nu}} \epsilon^{-\frac{2}{1+3\nu}}
\Big)\quad\andd\quad N_{\rm GLMO}(\epsilon) = O\Big( M^{\frac{4}{1+3\nu}}D^{\frac{4(1+\nu)}{1+3\nu}}\varepsilon^{-\frac{4}{1+3\nu}}\Big). %\qquad \mbox{for all}\;\; \nu\in[0,1]. 
\label{eq:comp_UCGS_nu}
\end{align}
In fact, it has been long known that the FO complexity in~\eqref{eq:comp_UCGS_nu} is optimal in the high-dimensional regime (cf.~\cite{Nemi_85}). In this respect, among the first-order AH methods, the universal CGS method achieves the optimal FO complexity  for obtaining  an $\epsilon$-optimal solution of~\eqref{eq:P}. More recently, the author in~\cite{Zhao_26a} proposed and analyzed a different first-order AH method with a single-loop structure,  and showed that for all $ \nu\in[0,1],$ the FO and GLMO complexities of this method are as follows: 
\begin{align}
N_{\rm FO} (\epsilon)= O\Big(M^{\frac{2}{1+\nu}}D^2\epsilon^{-\frac{2}{1+\nu}}\Big)\quad\andd\quad N_{\rm GLMO}(\epsilon) = O\Big(M^{\frac{2}{1+\nu}}D^2\epsilon^{-\frac{2}{1+\nu}}\Big). %\qquad \mbox{for all}\;\; \nu\in[0,1]. 
\label{eq:comp_PD_nu}
\end{align}
In view of~\eqref{eq:opt_GLMO_AH_nu}, it is clear that the method in~\cite{Zhao_26a} achieves the optimal GLMO complexity among the first-order AH methods (in the high-dimensional regime). 

At this moment, we know that in the high-dimensional regime, to obtain  an $\epsilon$-optimal solution of~\eqref{eq:P}, the optimal FO and GLMO complexties of the first-order AH methods  are
\begin{align}
O\Big( M^{\frac{2}{1+3\nu}} D^{\frac{2(1+\nu)}{1+3\nu}} \epsilon^{-\frac{2}{1+3\nu}}
\Big)\quad\andd\quad  O\Big(M^{\frac{2}{1+\nu}}D^2\epsilon^{-\frac{2}{1+\nu}}\Big), \quad \mbox{respectively}, \quad \mbox{for all }\, \nu\in[0,1]. 
\label{eq:opt_AH_nu}
\end{align}
Unfortunately,  neither method in~\cite{Ouyang_23} and~\cite{Zhao_26a} achieves the the optimal FO and GLMO complexties in~\eqref{eq:opt_AH_nu} {\em simultaneously} for all $\nu\in[0,1]$. That said, %we should 
it should be remarked that at the two ``endpoint'' cases, namely $\nu=0$ and $\nu=1$, there indeed exist several first-order AH methods that achieve the optimal FO and GLMO complexties simultaneously. Specifically, for $\nu=1$, the optimal FO and GLMO  complexties have been achieved by both the CGS  method~\cite{Lan_16} and the universal CGS method~\cite{Ouyang_23}; for $\nu=0$, the optimal FO and GLMO  complexties have been achieved by at least three different first-order AH methods, which are proposed in~\cite{Theku_20},~\cite{Asgari_22} and~\cite{Zhao_26a}, respectively. %That said, 
In contrast to these two ``endpoint'' cases, simultaneously achieving the the optimal FO and GLMO complexties in~\eqref{eq:opt_AH_nu} {\em for all $\nu\in [0,1]$}  is a considerably harder task. % uniformly}  is  a 
% has achieved 
To our knowledge, this has not been achieved by any first-order AH method in the literature. 

In this work, we fill this gap by developing a new first-order AH method for obtaining  an $\epsilon$-optimal solution of~\eqref{eq:P} that achieves the optimal FO and GLMO complexties in~\eqref{eq:opt_AH_nu} {\em simultaneously for all} $\nu\in[0,1]$. This method is based on our newly developed %the idea of 
{\em recentered homotopy smoothing framework}, which judiciously exploits various properties of the Moreau envelope and the curvature of the H\"older smooth $f$ (cf.~\eqref{eq:holder}). The rest of this paper is organized as follows. In Section~\ref{sec:prelim}, we provide some preliminaries of the H\"older smooth functions and Moreau envelopes. Next, in Section~\ref{sec:restart}, we analyze introduce the recentered homotopy smoothing framework and analyze its oracle complexities. We start by providing an overview of this framework and its subroutines (cf.\ Section~\ref{sec:overview}), and then analyze the iteration complexity of each subroutine (cf.\ Section~\ref{sec:triple_ave} to~\ref{sec:algo_iAPG}). Finally, we analyze the total LO and GLMO complexities of our framework (cf.\ Section~\ref{sec:total_comp}). 

%first-order primal-dual splitting method 
%However, the GLMO complexity does not match the optimal GLMO complexity of the class of AH methods as shown in~\eqref{eq:opt_comp_holder}. 

%In fact, this complexity is achieved by a single-loop AH method with with FO access of $f$. 

%Now, 
\begingroup
\renewcommand{\thealgorithm}{A}
\begin{algorithm}[t]
\caption{Active-Hull Methods}
\label{algo:LMO-AH}
\begin{algorithmic}
    \State {\bf Input}: $y_0 \in \dom g$, total number of iterations $K\ge 1$ and problem parameters (e.g., $M$, $\nu$, $D$)
    \For{\(k=0,1,2,\ldots,K-1\)}
%        \State Define the linear function \(\langle p_k,\cdot\rangle\).
        \State Construct dual variable $\lambda_k\in \bbR^n$ %and %call the LMO to 
        \State Compute $y_{k+1} \in \argmin_{y\in\bbR^n} \; \ip{\lambda_k}{y} + g(y)$
   \EndFor
        \State {\bf Output}: $\bary_K \in \conv\{y_0,\ldots,y_K\}$
\end{algorithmic}
\end{algorithm}
\endgroup

\section{Preliminaries} 
\label{sec:prelim}

%Define 
%\begin{equation}
%D_f(y,x):= f(y) - f(x) - \ipt{\nabla f(x)}{y-x}.
%\end{equation}

Let us first introduce the notion of $(\delta,L)$-inexact FO as proposed in~\cite{Dev_14}. As we shall see, this notion plays a critical role in the analyses of our algorithms.  

\begin{definition}[$(\delta,L)$-Inexact FO~\cite{Dev_14}] \label{def:delta_L_model}
%For any $u\in\bbR^n$, 
Let $h:\bbR^n\to \barbbR$ be a p.c.c.\ function. We say that $h$ is equipped with a $(\delta,L)$-inexact FO on $\bbR^n$ if at any $u\in\bbR^n$, the oracle outputs a pair $(\beta,l)\in\bbR\times\bbR^n$ such that % is called 
%first-order-approximation (FOA) 
%of $h$  if %for all $x\in\calU$, we have
\begin{equation}
0\le h(x) - \beta - \ipt{l}{x-u}\le %\beta + \ipt{l}{x-u} +
 (L/2)\normt{x-u}^2 + \delta, \quad \forall\,x\in\bbR^n. \label{eq:delta_L}
\end{equation}
%If we can compute such a $(\delta,\barL)$-FOA of $h$ at any $u\in\calU$, then  we say that $h$ is equipped with the $(\delta,\barL)$-inexact model on $\calU$.  %\Halmos
\end{definition}

Next, we record some standard results regarding the H\"older smooth function $f$. % in the following lemma. 
The proof of these results  can be found in~\cite[Proof of Lemma~1]{Nest_15} and~\cite[Proof of Definition~3.1]{Kerd_21}. 

\begin{lemma} \label{lem:holder_as_inexact}  %[{}]
For all $\nu\in[0,1)$ and any $\delta>0$, we have that for all $x,y\in\bbR^n$, 
\begin{equation}
f(y) - f(x) - \ipt{g_x}{y-x}\le \frac{L_\delta}{2}\normt{y-x}^2 + {\delta}, \quad \forall\,g_x\in \partial f(x),\quad \where\quad L_\delta:= \delta^{-\frac{1-\nu}{1+\nu}} M^{\frac{2}{1+\nu}}.  \label{eq:inexact_quad_ub}
 %\forall\,x,y\in\bbR^n, 
\end{equation}
%where $L_\delta:= \delta^{-\frac{1-\nu}{1+\nu}} M^{\frac{2}{1+\nu}}$. 
Consequently, we have 
\begin{equation}
 f(y)-  f(x) - \ip{g_x}{y-x}+ {\delta} \ge   \frac{1}{2L_\delta}\norm{g_y - g_x}^{2}  , \quad \forall\,g_x\in \partial f(x),\; \forall\,g_y\in \partial f(y), \quad \forall\, x,y\in \bbR^n. %\quad \where . 
\label{eq:inexact_lb}
\end{equation}
\end{lemma}
In view of Definition~\ref{def:delta_L_model}, the first part of Lemma~\ref{lem:holder_as_inexact} indicates that for any $x\in\bbR^n$ and any $g_x\in \partial f(x)$, $(f(x), g_x)$ is a $(\delta,L_\delta)$-inexact FO of $f$ at $x$, for any $\delta>0$.

Since our framework is based Moreau smoothing~\cite{Moreau_65}, we shall also introduce some basic properties of the Moreau envelope (see e.g.,~\cite[Section~12.4]{Bauschke_11}). To start with, let $\phi:\bbR^n\to \barbbR$ be a p.c.c.\ function and $\lambda>0$, and  define %, and for all $\lambda>0$, define 
\begin{alignat}{2}
\phi_\lambda(x) &:= \min_{u\in\bbR^n} \; \phi(u) + \frac{1}{2\lambda}\normt{u-x}^2,&\quad\forall\, x\in\bbR^n, \label{eq:phi_lambda}\\  %\;\andd\;
 \prox_{\lambda \phi}(x) &:= \argmin_{u\in\bbR^n} \; \phi(u) + \frac{1}{2\lambda}\normt{u-x}^2,&\quad \forall\, x\in\bbR^n. 
\end{alignat}
We call $\phi_\lambda$ the {\em Moreau envelope} of $\phi$ of parameter $\lambda$ and 
$\prox_{\lambda \phi}$  the {\em proximal operator} of $\lambda \phi$.

\begin{lemma} \label{lem:moreau}
%Let $\phi:\bbR^n\to \barbbR$ be a p.c.c.\ function and $\lambda>0$.  Then 
The function  $\phi_\lambda$ is convex and differentiable on $\bbR^n$ with 
\begin{equation}
\nabla \phi_\lambda(x) = \frac{x - \prox_{\lambda \phi}(x)}{\lambda},\quad  \forall\, x\in\bbR^n. 
\end{equation}
Furthermore, $\nabla \phi_\lambda$ is $\lambda^{-1}$-Lipschitz on $\bbR^n$.  In addition, 
for all $x\in \dom \partial\phi:=\{z\in\bbR^n: \partial\phi(z)\ne \emptyset\}$, 
we have 
\begin{equation}
0\le \phi(x) - \phi_\lambda(x) \le \frac{\lambda}{2}\normt{g_x}^2, \quad \forall\, g_x\in\partial \phi(x). 
%, \quad\forall\, x\in \dom \partial\phi. 
\label{eq:smoothing_bias}
\end{equation}
\end{lemma}

%\begin{lemma}
%For any $\nu\in(0,1]$, 
%%For all $x,y\in\bbR^n$, 
%we have 
%\begin{equation}
%D_f(y,x)\ge \frac{\nu}{1+\nu} M^{-\frac{1}{\nu}} \normt{g_y-\nabla f(x)}^{\frac{1+\nu}{\nu}}, \quad \forall\, x,y\in\bbR^n. 
%\end{equation} 
%\begin{equation}
%c_\nu
%\end{equation}
%\end{lemma}

%\paragraph{Notations.} 

\section{A Recentered Homotopy Smoothing Framework and Its Complexity Analysis}\label{sec:restart}

Our recentered homotopy smoothing framework is conceptually simple, %at a high level, 
but its implementation consists of several carefully designed subroutines. The complete set of algorithms %are which 
are shown in Algorithms~\ref{algo:restart} to~\ref{algo:FW}. In this section, we shall first provide an overview of these algorithms, followed by detailed description and complexity analysis of each one. Before proceeding, we remark that throughout this work, %we 
the problem parameters, including $M$, $\nu$ and $D$, are assumed to be known, and so we do not explicitly mention them as the input to any of the algorithms. %Developing a parameter-free version of of our method is certainly an important problem, and we leave it to future work. 
Also, for convenience, %or any $x\in\bbR^n$, 
define
\begin{equation}
\gap(x) := P(x) - P^*, \quad \forall\,x\in\bbR^n. 
\end{equation}
Also, in the sequel, let $x^*\in\calX^*$ be any optimal solution of~\eqref{eq:P}. 
%Given any p.c.c.\ function $P$ such that $P^*:=\min_{x\in\bbR^n}P(x) >-\infty $, define $\gap(x) := P(x) - P^*$ 

\subsection{Overview of Algorithms~\ref{algo:restart} to~\ref{algo:FW}} \label{sec:overview} %Description of Algorithms~\ref{algo:restart} and~\ref{algo:iAPG} }

As shown in Algorithm~\ref{algo:restart}, our recentered homotopy smoothing framework consists of $S\ge 1$ stages, ranging from $0$ to $S-1$. %namely $\{0,1,\ldots,S-1\}$. % ranging from $0$ to $S-1$. Ata 
First, note that the pre-start step in Algorithm~\ref{algo:restart} ensures that $\gap(a_0)\le\varepsilon_0$, as shown in the lemma below. 

\begin{lemma}
In Algorithm~\ref{algo:restart}, for any $x\in \dom g$, let $a_0$ and $\varepsilon_0$ be defined in the pre-start step. Then we have $\gap(a_0)\le\varepsilon_0$. 
%Fix any $x\in \dom g$ and pick any $g_x\in\partial f(x)$. Let $v\in \argmin_{u\in\bbR^n}\, \ip{g_x}{u} + g(u)$. Then for any $\delta>0$,  we have 
%\begin{equation}
%\gap(v)\le \frac{L_\delta}{2}\normt{y-x}^2 + {\delta},\qquad \mbox{where\;\; $L_\delta$ is defined in~\eqref{eq:inexact_quad_ub}.}\nn 
%\end{equation}
\end{lemma}

\begin{proof}
Let $G(x):= \ip{g_x}{x-a_0} + g(x) - g(a_0)$ be the (generalized) Wolfe-gap at $x$, and note that 
\begin{equation}
G(x)\ge \ip{g_x}{x-x^*} + g(x) - g(x^*) \ge f(x) - f(x^*) + g(x) - g(x^*) = P(x) - P^*. \label{eq:lb_G}
\end{equation}
%where $x^*\in\calX^*$. 
From~\eqref{eq:inexact_quad_ub}, for any $\delta>0$,  we have 
\begin{align}
P(a_0) = f(a_0) + g(a_0) &\le f(x) + \ipt{g_x}{a_0-x} +\frac{L_\delta}{2}\normt{a_0-x}^2 + {\delta} + g(a_0)\\
&=  f(x) + g(x) -G(x) +\frac{L_\delta}{2}\normt{a_0-x}^2 + {\delta} \\
&\le P^* +\frac{L_\delta}{2}\normt{a_0-x}^2 + {\delta},
\end{align}
where the last inequality follows from~\eqref{eq:lb_G}. This completes the proof. 
\end{proof}

Next, suppose that the start of some stage $0\le s\le S-1$, we have $a_s\in \dom g$ and $\varepsilon_s>0$ %at our disposal %, which satisfy  
that satisfy  $\gap(a_s)\le\varepsilon_s$. Based on $a_s$, we pick any $g_{s}\in \partial f(a_s)$ and define
\begin{equation}
\phi_s(u):= \phi(u) - \ip{g_s}{u}\quad\andd\quad \barg_s(u):= g(u) + \ipt{g_s}{u}, \quad \forall\,u\in\bbR^n. \label{eq:def_phi_s} 
\end{equation}
By choosing a suitable smoothing parameter $\lambda_s$,  let $\phi_{\lambda_s}$ be the {Moreau envelope} of $\phi$ of parameter $\lambda_s$ (cf.~\eqref{eq:phi_lambda}), and we invoke {\tt InexactAPG}, which is an inexact accelerated proximal gradient (APG) method (see e.g.,~\cite{Dev_14}) and shown in Algorithm~\ref{algo:iAPG}, to approximately solve the following (partially) smoothed problem: 
\begin{equation}
{\min}_{u\in\bbR^n}\;\; \phi_{\lambda_s}(u) + \barg_s(u) %\label{eq:prob_Phi_lambda}
\end{equation}  
The goal of %Algorithm~\ref{algo:iAPG} 
{\tt InexactAPG} is to return $(a_{s+1},\varepsilon_{s+1})$ such that %$\varepsilon_{s+1} $ and  
$\gap(a_{s+1})\le\varepsilon_{s+1}:= \varepsilon_s/2$. %As a result, 
If this is the case, then by induction, we know that  $\gap(a_s)\le \varepsilon_s= 2^{-s} \varepsilon_0$ for all $0\le s\le S-1$. % in Algorithm~\ref{algo:restart}. 

In the above, at each stage $s$, we call $a_s$ the ``center'', since 
\begin{equation}
0\in \partial  \phi_s(a_s)\quad\Longleftrightarrow\quad a_s\in {\argmin}_{u\in\bbR^n}\, \phi_s(u). 
\end{equation}
%$0\in \partial  \phi_s(a_s)$, or equivalently,  $a_s\in \argmin_{u\in\bbR^n} \phi_s(u)$. 
For this reason, the operation of defining $\phi_s$ and $\barg_s$ as ``shifted versions'' of $f$ and $g$ (cf.~\eqref{eq:def_phi_s}) 
 is called ``centering''. This operation may seem a bit unnatural at first glance, %on the surface level, 
 however, it plays a crucial role in bounding the (Moreau) smoothing error by exploiting % the curvature of 
the H\"older smoothness of $f$ (see Lemma~\ref{lem:smoothing_ub} %in Section~\ref{sec:algo_iAPG} 
below for details).  
%and $g_s$  

Next, let us turn our attention to Algorithm~\ref{algo:iAPG}. As aforementioned, given a pair $(a,\varepsilon)$ such that $\gap(a)\le \varepsilon$, Algorithm~\ref{algo:iAPG} runs an inexact APG method %to (approximately) solve
to solve the (partially) smoothed problem~\eqref{eq:prob_Phi_lambda}, with the goal to output a new pair $(a^+,\varepsilon^+)$ such that $\gap(a^+)\le \varepsilon^+:=\varepsilon/2$. Now, a careful reader may immediately spot the following gap: $\gap(a^+)$ is defined w.r.t.\ the original problem~\eqref{eq:P}, rather than~\eqref{eq:prob_Phi_lambda}, which replaces $\phi$ by $\phi_\lambda$. Naturally, in order to bound $\gap(a^+)$, we need to bound the smoothing error $\phi - \phi_\lambda$ at any $y\in\bbR^n$, which leads to the following lemma.

\begin{lemma} \label{lem:smoothing_ub}
Let $\phi$ and $\barg$ be defined in Line~\ref{line:def} in Algorithm~\ref{algo:iAPG}. %Then f
For any $\delta>0$, we have
\begin{equation}
\phi(y)  - \phi_\lambda(y)\le 2{\lambda}L_\delta (\gap(a) + \gap(y)+\delta), \quad \forall\, y\in\bbR^n. 
\end{equation}
\end{lemma}

\begin{proof}
Fix any $s_y\in \partial\phi(y)$. Since $\partial\phi = \partial f - g_a$, there exists $g_y\in\partial f(y)$ such that $s_y = g_y - g_a$. 
By~\eqref{eq:smoothing_bias}, %for any $s_y\in \partial\phi(y)$, 
we have 
\begin{equation}
\phi(y)  - \phi_\lambda(y)\le \frac{\lambda}{2} \normt{s_y}^2 = \frac{\lambda}{2} \normt{g_y - g_a}^2  \le \lambda\big(\normt{g_y - g_*}^2 + \normt{g_a - g_*}^2\big). \label{eq:phi-phi_lambda}
\end{equation}
By the first-order optimality condition, we know that there exists $g_*\in \partial f(x^*)$ such that $-g_*\in \partial g(x^*)$, and hence for any $\delta>0$, by~\eqref{eq:inexact_lb}, we have 
%$-g_*\in \partial g(x^*)$, and by~\eqref{eq:inexact_lb}, we have  
%and in addition,
\begin{align}
\normt{g_y - g_*}^2 &\le 2L_\delta(f(y)-  f(x^*) - \ip{g_*}{y-x^*}+ \delta/2)\\
&\le 2L_\delta(f(y)-  f(x^*) +g(y) - g(x^*)+ {\delta/2}) = 2L_\delta(P(y)-  P^*+ \delta/2). \label{eq:ub_gy_g*}
\end{align}
By the same reasoning, we have $\normt{g_a - g_*}^2\le 2L_\delta(P(a)-  P^*+ \delta/2)$. Combining this with~\eqref{eq:phi-phi_lambda} and~\eqref{eq:ub_gy_g*}, we complete the proof. 
\end{proof}

For convenience, define   
\begin{equation}
%\gap(x) := P(x) - P^*\qquad \andd\qquad 
\gap_\lambda(y):= \Phi_\lambda(y) - \Phi_\lambda^*, \quad\forall\, y\in\bbR^n. 
\end{equation}
%Based on this lemma, 
Note that we can easily decompose $\gap(y)$ as follows:
\begin{align}
\gap(y) = P(y) - P^* &= \phi(y)  - \phi_\lambda(y) + \phi_\lambda(y)+ \barg(y) - P^*\\
&\le \phi(y)  - \phi_\lambda(y) + \Phi_\lambda(y)  - \Phi_\lambda^* = \phi(y)  - \phi_\lambda(y) + \gap_\lambda(y),
%\le 2{\lambda}L_\delta (\gap(a) + \gap(y)+\delta)+ \gap_\lambda(y). %  - \Phi_\lambda^*. 
\end{align}
where the inequality follows from 
\begin{equation}
\Phi^\lambda_*={\min}_{u\in\bbR^n}\; \phi_\lambda(u) + \barg(u)\le  {\min}_{u\in\bbR^n}\; \phi(u) + \barg(u)= {\min}_{u\in\bbR^n}\; P(u) = P^*. 
\end{equation}
By Lemma~\ref{lem:smoothing_ub},   we immediately have the following lemma.

\begin{lemma}\label{lem:bound_gap}
Fix any $y\in\bbR^n$. 
For any $\delta>0$ and $0<\lambda<1/(2L_\delta)$, we have 
\begin{equation}
\gap(y)\le \frac{2\lambda L_\delta\big(\gap(a) + \delta\big) + \gap_\lambda(y)}{1-2\lambda L_\delta}.  
\end{equation}
In particular, if $\gap(a)\le \varepsilon$ for some $\varepsilon>0$, %$a\not\in\calX^*$, 
then choose $\delta=\varepsilon$ and $\lambda = c/(2L_\varepsilon)$ for some $c\in(0,1)$, and we have 
\begin{equation}
\gap(y)\le \frac{2c\, \varepsilon  + \gap_\lambda(y)}{1-c}.  \label{eq:bound_gap}
\end{equation}
\end{lemma}
Indeed, by choosing $y = a^+$,~\eqref{eq:bound_gap} in Lemma~\ref{lem:bound_gap} provides us the justification for Algorithm~\ref{algo:iAPG}: if $c$ is chosen to be sufficiently small, and the inexact APG method can find a point $a^+$ such that $\gap_\lambda(a^+)$ is sufficiently small, then we can ensure that $\gap(a^+)\le \varepsilon/2$.

Next, let us provide some implementation details of Algorithm~\ref{algo:iAPG}. At a high level, Algorithm~\ref{algo:iAPG} can be regarded as the APG method in~\cite[Algorithm~3]{Tseng_08} for solving~\eqref{eq:prob_Phi_lambda} with $(\delta,L)$-inexact FO of $\phi_\lambda$ (cf.~Definition~\ref{def:delta_L_model}) and inexact proximal oracle of $\barg$. %, with three sequences of iterates. 
These two inexact oracles were implemented by the subroutines {\tt InexactFO} (cf.~Algorithm~\ref{algo:triple-ave}) and {\tt InexactProx} (cf.~Algorithm~\ref{algo:FW}), respectively.  %The general structure of Algorithm~\ref{algo:iAPG}
 Similar to~\cite[Algorithm~3]{Tseng_08},  Algorithm~\ref{algo:iAPG} has three main sequences of iterates, namely $\{x_k\}_{0\le k\le K}$, $\{y_k\}_{0\le k\le K}$ and $\{z_k\}_{0\le k\le K}$. The other two sequences $\{\xi_k\}_{0\le k\le K}$ and $\{c_k\}_{0\le k\le K}$ %play auxiliary roles in 
facilitate defining the proximal subproblems to be solved by {\tt InexactProx}. 
The details of {\tt InexactFO}  and {\tt InexactProx} will be provided in Sections~\ref{sec:triple_ave} and~\ref{sec:FW}, respectively.  %Algorithms~\ref{al

Finally, we show that Algorithm~\ref{algo:restart} is indeed an AH method, namely an instance of Algorithm~\ref{algo:LMO-AH}. 
For $s\ge 0$, let $\calV_s$ denote all the points returned by the GLMO of $g$ before $a_s$ is formed. Clearly, from Line~\ref{line:prestart_restart} of Algorithm~\ref{algo:restart}, we have $\calV_0= \{a_0\}$. Now, at any stage $s\ge 0$, $a_s$ is passed to Algorithm~\ref{algo:iAPG}, where $x_0 = a_s$. From Line~\ref{line:x_+} of Algorithm~\ref{algo:iAPG}, we know that $x_K \in%\conv\{x_0,z_1,\ldots,z_{K-1}\}= 
\conv\{a_s,z_1,\ldots,z_{K-1}\}$. Let $\calA_{s,k}$ denote the set of points returned by the GLMO of $g$ in {\tt InexactProx} (i.e.,~Algorithm~\ref{algo:FW}) at the $k$-th iteration of Algorithm~\ref{algo:iAPG}. Since $a_{s+1} =x_{K}$, we have $\calV_{s+1}= \big(\bigcup_{k=0}^{K-1}\calA_{s,k} \big)\bigcup\calV_s$. Since $\eta_0=1$ in Algorithm~\ref{algo:FW}, we know that $z_k\in \conv \calA_{s,k}$ for all $0\le k\le K-1$. As a result, if $a_s\in \conv\calV_s$, then $a_{s+1} = x_K\in \conv\{a_s,z_1,\ldots,z_{K-1}\} \subseteq \conv \calV_{s+1}$. Since apparently $a_0\in \conv\calV_0$, using induction, we know that $a_s\in \conv\calV_s$ for all $0\le s\le S$. This shows that Algorithm~\ref{algo:restart} is an AH method.  %Consequently, 

%
%\begin{prop}
%Algorithm~\ref{algo:restart} is an AH method. 
%\end{prop}
%
%\begin{proof}
%
%\end{proof}

%\begin{proof}
%First, note that from~\eqref{eq:smoothing_bias}, we know that  %for all $y\in\bbR^n$, 
%
%Define $\gap_\lambda:= \gap_{\Phi_\lambda}$. 
%%This, together with 
%
%
%%Also,  and hence %we have %note that
%%\begin{equation}
%%\normt{s_y}^2 = \normt{g_y - g_a}^2  \le 2\big(\normt{g_y - g_*}^2 + \normt{g_a - g_*}^2\big).
%%\end{equation}
%Since $0<\lambda<1/(2L_\delta)$, by rearranging this inequality,  we complete the proof. 
%\end{proof}

\setcounter{algorithm}{0}
\begin{algorithm}[t]
\caption{Recentered Homotopy Smoothing Framework}
\label{algo:restart}
\begin{algorithmic}[1]
    \State {\bf Input}: $x_\inp\in\dom g$, total number of stages $S\ge 1$
    \State \label{line:prestart_restart} {\bf Pre-start}:  Pick any $g_x\in\partial f(x)$ and  let $a_0\in \argmin_{u\in\bbR^n}\, \ip{g_x}{u} + g(u)$. Fix any $\delta>0$ and let 
\begin{equation}\setlength{\abovedisplayskip}{4pt}\setlength{\belowdisplayskip}{4pt}
\varepsilon_0:=({L_\delta}/{2})\normt{a_0-x}^2 + {\delta},\qquad \mbox{where $L_\delta$ is defined in~\eqref{eq:inexact_quad_ub}.}\nn 
\end{equation}    
%    $({L_\delta}/{2})\normt{y-x}^2 + {\delta}$ %$a_0 \in \dom g$ and $\varepsilon_0>0$ such that $\gap(a_0)\le \varepsilon_0$  %and 
    %\State {\bf Choose}: $\{\varepsilon_k\}_{k\ge 0}$ such that $\gap(a_0)\le \varepsilon_0$ %, $\{\lambda_k\}_{k\ge 0}$
     \For{\(s=0,1,2,\ldots,S-1\)}
        \State %($a_{k+1}$, $\varepsilon_{k+1}$) 
        $(a_{s+1},\varepsilon_{s+1}) = {\tt InexactAPG}(f,g, a_s, \varepsilon_s)$ %, \varepsilon_{k+1})$ %, \lambda_k)$
   \EndFor
        \State {\bf Output}: $a_S $
\end{algorithmic}
\end{algorithm}

\begin{algorithm}[t]
\caption{Inexact %Accelerated Proximal Gradient 
APG Method ${\tt InexactAPG}(f,g,a, \varepsilon)$}
\label{algo:iAPG}
\begin{algorithmic}[1]
    \State {\bf Require}: $%P(a) - P^*
    \gap(a)\le \varepsilon$ %, $0<\varepsilon^+< \varepsilon$
        %        $\lambda = 1/(18L_\varepsilon)$, $\bar\delta=(2/9)\varepsilon$ %, $\varepsilon^+= \varepsilon/2$, 
          \State {\bf Choose}: Total number of iterations $K\ge 1$, smoothing parameter $\lambda>0$,  weight sequence $\{A_k\}_{0\le k\le K}$, inner accuracy parameters $\delta>0$ and        $\{\tau_k\}_{0\le k\le K}$. % and such that 
    \State {\bf Definition}: \label{line:def}
    Fix any $g_a\in \partial f(a)$, and let $\phi:= f-\ipt{g_a}{\cdot}$ and  $\barg:= g + \ipt{g_a}{\cdot}$. Also, let 
\begin{equation}
\textstyle \alpha_k := A_{k+1} - A_k %:=  \sum_{i=0}^{k-1} \alpha_i
\quad \andd\quad \theta_k := \alpha_k/A_{k+1},\quad \forall\, 0\le k\le K-1.    \label{eq:def_alpha_theta}
\end{equation}
%\State {\bf Set}: 
      
     %$A_k:= \sum_{i=0}^{k-1} \alpha_i$  and $\theta_k := \alpha_k/A_{k+1}$.  
    %\vspace{1ex}
    \Statex %{\bf Overview}: 
        Run the following inexact APG method to solve 
\begin{equation}
\Phi_\lambda^*:={\min}_{u\in\bbR^n}\;\; \big\{\Phi_\lambda(u):= \phi_\lambda(u) + \barg(u)\big\} \label{eq:prob_Phi_lambda}
\end{equation}        
%        $$:
%        \vspace{1ex}
    \State {\bf Initialize}: $x_0:=z_0 :=a$, $\xi_0:= 0$ %, $A_0 = 0$
     \For{\(k=0,1,2,\ldots,K-1\)}
        \State $y_k := (1-\theta_k) x_k + \theta_k z_k$
        \State $(w_k,l_k) := {\tt InexactFO}(\phi,a,\lambda,y_k,  \delta)$
        \State $\xi_{k+1}:= \xi_k + \alpha_k l_k$ %, 
        \State $c_{k+1} := a - \xi_{k+1}$
        \State $z_{k+1} := {\tt InexactProx}(\barg,c_{k+1}, A_{k+1},\tau_{k+1})$
        \State \label{line:x_+} $x_{k+1} := (1-\theta_k) x_k + \theta_k z_{k+1}$
       % \State $U_{k+1}: = \big(\sum_{i=0}^{k+1} A_i\tau_i + 2\delta\sum_{i=0}^{k+1} A_i + D^2/2\big)/A_{k+1}$ %+ \tau_{k+1} $
        %\If{$U_{k+1}\le \bar\delta$}
        %\EndIf
   \EndFor
   \State \Return $(x_{K},\varepsilon/2)$
        %\State {\bf Output}: $\bary_K \in \conv\{y_0,\ldots,y_K\}$
\end{algorithmic}
\end{algorithm}

\begin{algorithm}[t]
\caption{Inexact First-Order Oracle %Accelerated Gradient Method 
$%(\beta,l) := 
{\tt InexactFO}(\phi,a,\lambda,y, \delta)$}
\label{algo:triple-ave}
\begin{algorithmic}[1]
\State {\bf Require}: $0\in \partial \phi(a)$
\vspace{1ex}
  \State \label{line:nu=1}  If $\nu=1$, then return $%(\beta,l) := 
  (\phi(y),\nabla \phi(y))$; otherwise, run the following %accelerated %dual averaging 
  ``triple-averaging''   method to approximately solve 
\begin{equation}
%\phi_\lambda(y):=
Q_\lambda^*:={\min}_{u\in\bbR^n}\; \big\{Q_\lambda(u):= \phi(u) + \omega_\lambda(u)\big\}, \quad\where\quad  \omega_\lambda(u):= (2\lambda)^{-1}\normt{u-y}^2.  \label{eq:Q_lambda}
\end{equation}  

   %\State {\bf Choose}: $\{\vartheta_j\}_{j\ge 0}$
%    \State {\bf Define}:  $B_0 := 1$, $B_j= B_0+\sum_{i=0}^{j-1} b_i $, $\vartheta_j := b_j/B_{j+1}$

    \State {\bf Initialize}: $p_0:=a$, $r_0 :=y$, $(\beta_0,l_0) := (\phi(a),0)$, $\psi_0(u) := \phi(a) + \omega_\lambda(u)$,     $\psi_0^* :=\psi_0(r_0)=\phi(a)$, 
     $\Delta_{0}:=Q_\lambda(p_0)- \psi_0^*= (2\lambda)^{-1}\normt{a-y}^2$ % Q_\lambda(p_{0}) -  \psi^*_{0}$
    %$\psi_0(u) := (1/2) \normt{u-y}^2$ %, $A_0 = 0$
     \For{\(k=0,1,2,\ldots\)}
     \If{$\Delta_{j}\le \delta$}
          \State \Return $(\psi^*_{j},l_{j})$
        \EndIf
     \State \label{line:vartheta} $\vartheta_j = \min\{1/2, (\zeta_j/2)^{1/w}\}$, where $\zeta_j := \lambda^{-1} M^{-\frac{2}{1+\nu}}(\Delta_j/2)^{\frac{1-\nu}{1+\nu}}  $ and $w := \frac{1+3\nu}{1+\nu}$
        \State \label{line:q} $q_j := (1-\vartheta_j) p_j + \vartheta_j r_j$
        \State Pick any $s_j\in \partial \phi(q_j)$
         \State \label{line:beta} $\beta_{j+1} := (1-\vartheta_j) \beta_j + \vartheta_j (\phi(q_j) - \ipt{s_j}{q_j})$
        \State \label{line:l} $l_{j+1} := (1-\vartheta_j) l_j + \vartheta_j s_j$
%        \State $m_j(u):= \phi(q_j) + \ipt{\nabla \phi(q_j)}{u-q_j} + (2\lambda)^{-1}\normt{u-y}^2$
%        \State $\psi_{j+1}(u) := \beta_{j+1} + \ipt{l_{j+1}}{u} + (2\lambda)^{-1}\normt{u-y}^2$ 
        \State \label{line:r} $r_{j+1} := \argmin_{u\in\bbR^n}\, \{ \psi_{j+1} (u) := \beta_{j+1} + \ipt{l_{j+1}}{u} + \omega_\lambda(u)\}$
        \State \label{line:psi*} $\psi^*_{j+1} := \psi_{j+1} (r_{j+1})$
        \State \label{line:p} $p_{j+1} := (1-\vartheta_j) p_j + \vartheta_j r_{j+1}$
		\State $\Delta_{j+1}: = Q_\lambda(p_{j+1}) -  \psi^*_{j+1}$
        
   \EndFor
        %\State {\bf Output}: $\bary_j \in \conv\{y_0,\ldots,y_K\}$
\end{algorithmic}
\end{algorithm}

\begin{algorithm}[t]
\caption{Inexact Proximal Oracle ${\tt InexactProx}(\barg,c, A,\tau)$}
\label{algo:FW}
\begin{algorithmic}[1]
  \Statex  Run the following generalized Frank-Wolfe method to approximately solve 
\begin{equation}
\barQ_*:={\min}_{u\in\bbR^n}\; \big\{\barQ(u):= Q(u) +\barg(u) \big\},\quad \where\quad  %Q:= (2A)^{-1}\normt{\cdot-c}^2. 
Q(u):= (2A)^{-1}\normt{u-c}^2,\;\;\; \forall\,u\in\bbR^n. \label{eq:FW_subprob}
\end{equation}  

    \State {\bf Initialize}: $u_0\in \dom\barg$ %, $\psi_0(u) \equiv 0$, $d_0:= Q_\lambda(q_0)$ 
    %$\psi_0(u) := (1/2) \normt{u-y}^2$ %, $A_0 = 0$
     \For{\(j=0,1,2,\ldots\)}
       \State \label{line:def_v_j} Compute %$g_j:=\nabla f(x_j)$, 
       $v_{j} \in \argmin_{u\in \bbR^n}\{\ip{\nabla Q(u_j)}{u} + \barg(u)\}$ %and  $\eta_j \in [0,1]$
		\State Compute $G_j:=\ip{\nabla Q(u_j)}{u_j-v_j} + \barg(u_j) - \barg(v_j)$
		\If{$G_j\le \tau$}
		\State \Return $u_j$
		\EndIf
		%\Else 
		\State Compute $u_{j+1} := (1-\eta_j) u_j + \eta_j v_j$, where   $\eta_j =2/(j+2)$. 
		%\EndIf
   \EndFor
        %\State {\bf Output}: $\bary_j \in \conv\{y_0,\ldots,y_K\}$
\end{algorithmic}
\end{algorithm}

\subsection{Description and Analysis of Algorithm~\ref{algo:triple-ave}} \label{sec:triple_ave}

The goal of Algorithm~\ref{algo:triple-ave} is to provide a $(\delta,L)$-inexact FO of $\phi_\lambda$ at $y$, and this is achieved separately in two cases. The first case is $\nu=1$, where $(\phi(y),\nabla \phi(y))$ is returned directly (cf.~Line~\ref{line:nu=1}); the second case is $0\le \nu<1$, where the ``triple-averaging'' (TA)  method in~\cite[Algorithm~5]{Burns_26} is used to solve the Moreau smoothing subproblem %as shown 
in~\eqref{eq:Q_lambda}. %Despite the same method, 
It is %important to mention 
worth mentioning that the authors in~\cite{Burns_26} analyzed the TA method for solving~\eqref{eq:Q_lambda} under a very different setting than ours. Specifically, in their setting, $\phi_\lambda$ is $M$-smooth on $\dom\barg$ (i.e., $\nu=1$) and  % our problem~\eqref{eq:Q_lambda} is very different from the one , where %--- in that problem 
%the TA method was analyzed under a very different setting in~\cite{Burns_26} than ours: specifically, 
%$\nu=1$ and 
$\omega_\lambda$ is a general $\lambda^{-1}$-strongly-convex function (i.e., not necessarily equal to $(2\lambda)^{-1}\normt{\cdot-y}^2$).

In our implementation, %The algorithm 
we maintain a ``gap'' sequence $\{\Delta_j\}_{j\ge 0}$, and the TA method is terminated  as soon as $\Delta_j\le \delta$. (Indeed, our notion of ``gap'' is akin to the ``primal-dual certificate'' in~\cite[Definition~2.4]{Burns_26}.)  This termination criterion is justified in Lemma~\ref{lem:delta_L_inexact} below.

\begin{lemma}\label{lem:delta_L_inexact}
Let $\nu\in[0,1)$. 
Upon termination of Algorithm~\ref{algo:triple-ave}, we obtain $(\psi^*,l)\in\bbR\times \bbR^n$ that serves as a $(2\delta,2/\lambda)$-inexact FO of $\phi_\lambda$ at $y$, namely %  such that 
\begin{align}
0\le \phi_\lambda(x)  - \psi^* - \ip{l}{x-y} \le \lambda^{-1} \normt{x-y}^2 + 2\delta,  \quad \forall\, x\in\bbR^n. 
\end{align}
%where $\epsilon = \delta$ if $\nu\in[0,1)$ and $\epsilon = 0$ if $\nu=1$. 
\end{lemma}

% but before introducing it, 
The proof of Lemma~\ref{lem:delta_L_inexact} leverages several results below. First, 
 let us provide  % which follows from from 
a general result that may be of independent interest.  
% when $0\le \nu<1$, with the goal to . (As shown in Line~\ref{line:nu=1}, if 

\begin{prop} \label{prop:inexact_model_moreau_env}
Let $\varphi:\bbR^n\to \barbbR$ be a p.c.c.\ function, and $(\beta,l)\in\bbR\times \bbR^n$ satisfy that 
\begin{equation}
\ell(u):=\beta + \ip{l}{u}\le \varphi(u), \quad \forall\, u\in\bbR^n.  \label{eq:ell_varphi}
\end{equation}
%$$ for all $u\in\bbR^n$. 
Define 
\begin{alignat}{2}
\varphi_\lambda(x) &:= {\min}_{u\in\bbR^n} \; \big\{%R_\lambda(u;x):= 
\varphi(u) + (2\lambda)^{-1}\normt{u-x}^2\}, \qquad &\forall\, x\in\bbR^n,\\
\ell_\lambda(x)&:= {\min}_{u\in\bbR^n}\; \big\{%\psi_\lambda(u):=
\ell(u) + (2\lambda)^{-1}\normt{u-x}^2\big\},\qquad &\forall\, x\in\bbR^n. 
\end{alignat}
%$\varphi_\lambda(x) := \min_{u\in\bbR^n} \; \varphi(u) + (2\lambda)^{-1}\normt{u-x}^2$ for $x\in\bbR^n$. 
%For any $y\in\bbR^n$, define 
%\begin{equation}
%%\psi^*:=
%\ell_\lambda(y):= {\min}_{u\in\bbR^n}\; \big\{\psi(u):= \beta + \ip{l}{u} + (2\lambda)^{-1}\normt{u-y}^2\big\}\quad\andd\quad \Delta = Q_\lambda(y)  -  \psi^*.
%\end{equation}
If $\varphi_\lambda(y)  -  \ell_\lambda(y)\le \delta$ for some $y\in \bbR^n$, then we have 
\begin{align}
0\le \varphi_\lambda(x)  - \ell_\lambda(y) - \ip{l}{x-y} \le \lambda^{-1} \normt{x-y}^2 + 2\delta,  \quad \forall\, x\in\bbR^n. 
\end{align}
\end{prop}

\begin{proof}\renewcommand{\qed}{}
%Fix any $x\in\bbR^n\renewcommand{\qed}{}.$
First, note that
\begin{equation}
\ell_\lambda(x) = \beta + \ip{l}{x} - (\lambda/2)\normt{l}^2 = \ell(x) - (\lambda/2)\normt{l}^2, \quad \forall\, x\in\bbR^n.
\end{equation}
As a result, by~\eqref{eq:ell_varphi}, we have 
\begin{equation}
\varphi_\lambda(x) \ge \ell_\lambda(x) = \ell_\lambda(y) + \ell_\lambda(x) - \ell_\lambda(y) = \ell_\lambda(y) + \ell(x) - \ell(y) = \ell_\lambda(y) + \ip{l}{x-y},\quad \forall\, x\in\bbR^n. 
\end{equation}
This, together with Lemma~\ref{lem:moreau}, implies that 
\begin{align}
0&\le \varphi_\lambda(x)  - \ell_\lambda(y) - \ip{l}{x-y} \\
&\le \varphi_\lambda(y)- \ell_\lambda(y) + \ip{\nabla \varphi_\lambda(y)-l}{x -y} + (2\lambda)^{-1}  \normt{x-y}^2,\label{eq:varphi_1}\\
&\le \varphi_\lambda(y)- \ell_\lambda(y) + (\lambda/2)\normt{\nabla \varphi_\lambda(y)-l}^2 + \lambda^{-1}  \normt{x-y}^2, \qquad \forall\, x\in\bbR^n, \label{eq:varphi_2}
\end{align}
where in~\eqref{eq:varphi_2} we use Young's inequality. 
Note that by minimizing~\eqref{eq:varphi_1} over $x\in\bbR^n$, we have % implies that 
$$(\lambda/2) \normt{\nabla \varphi_\lambda(y)-l}^2 \le \varphi_\lambda(y)- \ell_\lambda(y).$$ This, combined with~\eqref{eq:varphi_2} and that $\varphi_\lambda(y)  -  \ell_\lambda(y)\le \delta$, yields  
\begin{align}
\varphi_\lambda(x)  - \ell_\lambda(y) - \ip{l}{x-y}\le 2(\varphi_\lambda(y)- \ell_\lambda(y)) + \lambda^{-1}  \normt{x-y}^2\le 2\delta + \lambda^{-1}  \normt{x-y}^2. \tag*{$\square$} 
\end{align}
\end{proof}
\noindent 
In view of Definition~\ref{def:delta_L_model}, Proposition~\ref{prop:inexact_model_moreau_env} indicates that if $\varphi_\lambda(y)  -  \ell_\lambda(y)\le \delta$  for some $y\in \bbR^n$, then $(\ell_\lambda(y),l)$ is a $(2\delta,2/\lambda)$-inexact FO of $\varphi_\lambda$ at $y$. 
%Now, %for analysis purposes,
%to introduce Corollary~\ref{cor:delta_L_inexact}, we also need to define  
Next, we define 
%consider the following sequence of functions: 
\begin{equation}
\psi_{0}(u) := \beta_{0} + \ipt{l_{0}}{u} + \omega_\lambda(u),\quad\forall\,u\in\bbR^n.
\end{equation}
Also note that by Line~\ref{line:r}, % in Algorithm~\ref{algo:triple-ave}, 
we have %that  for all $j\ge 0$,
\begin{equation}
\psi_{j}(u) = \beta_{j} + \ipt{l_{j}}{u} + \omega_\lambda(u),\quad\forall\,u\in\bbR^n, \quad \forall\, j\ge 0.  \label{eq:def_psi_j}  
\end{equation}

The next lemma introduces some basic properties of the sequence of functions $\{\psi_j\}_{j\ge 0}$: 

\begin{lemma}\label{lem:psi_j}
In Algorithm~\ref{algo:triple-ave}, %for any $j\ge 0$, 
 we have 
\begin{align}
\begin{split}
&\psi_{j+1}(u) = (1-\vartheta_j)\psi_{j}(u) + \vartheta_j m_j(u),\quad  \\
&\qquad\where\quad  m_j(u):=\phi(q_j) + \ipt{s_j}{u-q_j} + \omega_\lambda(u)\le Q_\lambda(u),\quad \forall\,u\in \bbR^n,\quad\forall\,j\ge 0. 
\end{split}
%\andd
%&\psi_{j}(u) = B_j^{-1}\textstyle \sum_{i=0}^{j-1} \, b_i\big(\phi(q_i) + \ipt{\nabla\phi(q_i)}{u-q_i}\big) + \omega_\lambda(u), \quad && \forall\,u\in \bbR^n, 
\label{eq:recur_psi}
\end{align}
In addition, %for any $j\ge 0$,  
we have 
\begin{equation} 
%\quad \Longrightarrow\quad 
\psi_j(u) \le Q_\lambda(u), \quad\forall\,u\in \bbR^n, \qquad\andd\qquad 
r_{j} = {\argmin}_{u\in\bbR^n}\;\psi_{j}(u), \quad \forall\, j\ge 0.  \label{eq:def_r_j} 
\end{equation}
\end{lemma}

\begin{proof}
Note that for any $j\ge 0$, from Lines~\ref{line:beta} and~\ref{line:l}, we have 
\begin{align}
\psi_{j+1}(u) &= \beta_{j+1} + \ipt{l_{j+1}}{u} + \omega_\lambda(u)\\
&= (1-\vartheta_j) \beta_j + \vartheta_j \big(\phi(q_j) - \ipt{s_j}{q_j}\big) + (1-\vartheta_j) \ipt{l_j}{u} + \vartheta_j \ipt{s_j}{u} + \omega_\lambda(u)\\
&= \vartheta_j \big(\phi(q_j) + \ipt{s_j}{u-q_j} + \omega_\lambda(u)\big) + (1-\vartheta_j)\big(\beta_j + \ipt{l_j}{u}+ \omega_\lambda(u)\big)\\
&= \vartheta_j m_j(u) + (1-\vartheta_j)\psi_{j}(u). 
\end{align}
The fact that $m_j\le Q_\lambda$ on $\bbR^n$ follows from % the convexity of $\phi$ and that
 $s_j\in \partial \phi(q_j)$. Next, note that since $0\in\partial\phi(a)$, % and $\phi$ is convex on $\bbR^n$, 
we have $\phi(u)\ge \phi(a)$ for all $u\in\bbR^n$, and since $(\beta_0,l_0) := (\phi(a),0)$, we have %hence %$ $
\begin{equation}
\psi_0(u) = \phi(a) + \omega_\lambda(u)\le \phi(u) + \omega_\lambda(u) = Q_\lambda(u), \quad \forall\,u\in\bbR^n. 
\end{equation}
Based on~\eqref{eq:recur_psi}, we then use induction to conclude that $\psi_j\le Q_\lambda$ on $\bbR^n$ for all $j\ge 0$. Finally, note that for $j\ge 1$, that $r_{j} = {\argmin}_{u\in\bbR^n}\;\psi_{j}(u)$ follows from Line~\ref{line:r}. For $j=0$, since $r_0 = y$, %and $(\beta_0,l_0) := (\phi(a),0)$, 
we have 
\begin{equation}
{\argmin}_{u\in\bbR^n}\;\psi_{0}(u) = {\argmin}_{u\in\bbR^n}\;\phi(a) + \omega_\lambda(u) = y = r_0. \tag*{$\square$}
\end{equation}
\renewcommand{\qed}{}
%This completes the proof. 
%Suppose that~\eqref{eq:def_psi} holds for some $j\ge 0$. Then we have
%\begin{align}
%\psi_{j+1}(u) &= \beta_{j+1} + \ipt{l_{j+1}}{u} + \omega_\lambda(u)\\
%&= (1-\vartheta_j) \beta_j + \vartheta_j \big(\phi(q_j) - \ipt{\nabla\phi(q_j)}{q_j}\big) + (1-\vartheta_j) \ipt{l_j}{u} + \vartheta_j \ipt{\nabla\phi(q_j)}{u} + \omega_\lambda(u)\\
%&= \vartheta_j \big(\phi(q_j) + \ipt{\nabla\phi(q_j)}{u-q_j} + \omega_\lambda(u)\big) + (1-\vartheta_j)\big(\beta_j + \ipt{l_j}{u}+ \omega_\lambda(u)\big)\\
%&= B_{j+1}^{-1}\big\{b_j\big(\phi(q_j) + \ipt{\nabla\phi(q_j)}{u-q_j} + \omega_\lambda(u)\big)\\
%&\hspace{4cm} + \textstyle \sum_{i=0}^{j-1} \, b_i\big(\phi(q_i) + \ipt{\nabla\phi(q_i)}{u-q_i}\big) + B_j\omega_\lambda(u)\big\}\\
%&= B_{j+1}^{-1} \textstyle \sum_{i=0}^{j} \, b_i\big(\phi(q_i) + \ipt{\nabla\phi(q_i)}{u-q_i}\big) + \omega_\lambda(u). %\big \} 
%\end{align}
\end{proof}

Based on Proposition~\ref{prop:inexact_model_moreau_env} and Lemma~\ref{lem:psi_j}, we are ready to prove Lemma~\ref{lem:delta_L_inexact}.

\begin{proof}[Proof of Lemma~\ref{lem:delta_L_inexact}]
Fix any $j\ge 0$. Define $\ell^j(u):= \beta_{j} + \ipt{l_{j}}{u}$ for $u\in\bbR^n$, so that $\psi_j = \ell^j + \omega_\lambda$ on $\bbR^n$. From~\eqref{eq:def_r_j}, we know that $\ell^j(u)\le \phi(u)$ for $u\in\bbR^n$. 
Note that %for all $j\ge 0$, 
$Q_\lambda^* = \phi_\lambda(y)$, $\psi^*_{j} = \ell^j_\lambda(y)$ and so $\Delta_j = Q_\lambda(p_{j}) -  \psi^*_{j}\ge Q_\lambda^* - \psi^*_{j}= \phi_\lambda(y) -\ell^j_\lambda(y)$. Since Algorithm~\ref{algo:triple-ave} terminates as soon as $\Delta_j\le \delta$, % and $(\psi^*_{j},l_{j})$ is returned, 
we have $\phi_\lambda(y) -\ell^j_\lambda(y)\le \delta$ upon termination. Invoking Proposition~\ref{prop:inexact_model_moreau_env}, we complete the proof. 
\end{proof}

The next lemma justifies the choice of the stepsize $\vartheta_j$ in Line~\ref{line:vartheta}. 
\begin{lemma} \label{lem:lb_L_j}
In Algorithm~\ref{algo:triple-ave}, for any $j\ge 0$, define 
\begin{equation}
\delta_j:= \vartheta_j \Delta_j/2 \quad \andd\quad L_j:= L_{\delta_j} = \delta_j^{-\frac{1-\nu}{1+\nu}} M^{\frac{2}{1+\nu}}. \label{eq:def_delta_L}
\end{equation}
Then we have 
\begin{equation}
\frac{1-\vartheta_j}{\lambda \vartheta_j^2}\ge L_j. 
\end{equation}
\end{lemma}

\begin{proof}
Let us first show that by the choice of $\vartheta_j$ in Line~\ref{line:vartheta}, we have 
\begin{equation}
\frac{\vartheta_j^w}{1-\vartheta_j} \le %\lambda^{-1} M^{-\frac{2}{1+\nu}}(\Delta_j/2)^{\frac{1-\nu}{1+\nu}}
  \zeta_j. 
\end{equation}
Indeed, if $\vartheta_j = 1/2$,  then we have $(1/2)^w\le \zeta_j/2$, and hence ${\vartheta_j^w}/({1-\vartheta_j})= 2{(1/2)^w}\le \zeta_j$; otherwise, $\vartheta_j = (\zeta_j/2)^{1/w}$, and we have  $(\zeta_j/2)^{1/w}\le 1/2$, which then implies that 
%and hence 
\begin{equation}
\frac{\vartheta_j^w}{1-\vartheta_j} = \frac{\zeta_j/2}{1-(\zeta_j/2)^{1/w}}\le \zeta_j. 
\end{equation}
%${\vartheta_j^w}/({1-\vartheta_j}) = (\zeta_j/2)/(1-(\zeta_j/2)^{1/w}) $
Now, from the definition of $\zeta_j$ and $w$ in Line~\ref{line:vartheta}, we have 
\begin{align}
\frac{\lambda\vartheta_j^{2}}{1-\vartheta_j} =\frac{\lambda\vartheta_j^{w+\frac{1-\nu}{1+\nu}}}{1-\vartheta_j}\le \lambda\vartheta_j^{\frac{1-\nu}{1+\nu}} \zeta_j \le  M^{-\frac{2}{1+\nu}}(\vartheta_j\Delta_j/2)^{\frac{1-\nu}{1+\nu}} = M^{-\frac{2}{1+\nu}}\delta_j^{\frac{1-\nu}{1+\nu}} = L_j^{-1}. \tag*{$\square$}
\end{align}\renewcommand{\qed}{}
%This completes the proof. 
\end{proof}

Next, let us analyze the total number of iterations of Algorithm~\ref{algo:triple-ave}. 
For any $\delta>0$, define
\begin{equation}
J_\TA(\delta) = \min\{j\ge 0: \Delta_{j}\le \delta\}.
\end{equation}

%For convenience, for any $j\ge 0$, define

\begin{theorem} \label{thm:TA}
Let $\nu\in[0,1)$. 
For $j\ge 0$, we have
\begin{align}
%\frac{1-\vartheta_j}{\lambda\vartheta_j^2}\ge L_j = L_{\delta_j} = 
\Delta_{j+1}\le  \left(1 - ({1}/{2})\min\left\{{1}/{2},\, \left({\zeta_j}/{2}\right)^{1/w}\right\}\right)\Delta_{j},
\end{align}
where $\zeta_j$ and $w$ are defined in Line~\ref{line:vartheta}. 
Define $j_0:= \min\{0\le j \le J_\TA(\delta): \Delta_j\le H\}$, where $H:= 2^{\frac{1-3\nu}{1-\nu}} \lambda^{\frac{1+\nu}{1-\nu}}M^{\frac{2}{1-\nu}}$. Then we have
\begin{alignat}{2}
\Delta_{j} &\le (3/4)^j \Delta_0, &\qquad&\forall\, 0\le j\le  j_0,\\
\Delta_j&\le \big(H^{-\gamma} + \gamma c(j-j_0)\big)^{-1/\gamma},&\qquad &\forall\, j_0\le j\le  J_\TA(\delta),
\end{alignat}
where 
\begin{equation}
c:=2^{-\frac{3(1+\nu)}{1+3\nu}} \lambda^{-\frac{1+\nu}{1+3\nu}}M^{-\frac{2}{1+3\nu}}\quad \andd\quad \gamma:=\frac{1-\nu}{1+3\nu}. \label{eq:def_c_gamma}
\end{equation}
%$$ and $$. 
\end{theorem}

\begin{proof}
Fix any $j\ge 0$. 
Note that from Lines~\ref{line:q} and~\ref{line:p}, we have  
\begin{equation}
p_{j+1} - q_j = \vartheta_j(r_{j+1} - r_j). \label{eq:triangle_identity}
\end{equation}
Consequently, we have 
\begin{align}
\psi_{j+1}^* &= \psi_{j+1} (r_{j+1})\\
 &= (1-\vartheta_j)\psi_{j}(r_{j+1}) + \vartheta_j m_j(r_{j+1})\\
&\ge  (1-\vartheta_j)\big(\psi_{j}^*+{(2\lambda)}^{-1}\normt{r_{j+1} - r_j}^2\big) + \vartheta_j m_j(r_{j+1}) \label{eq:analysis_TA_1}\\
&\ge (1-\vartheta_j)\big(\psi_{j}^*+(2\lambda\vartheta_j^2)^{-1}\normt{p_{j+1} - q_j}^2\big) + \vartheta_j m_j(r_{j+1}) + (1-\vartheta_j)(m_j(p_{j}) - Q_\lambda(p_{j}))\label{eq:analysis_TA_2}\\
&\ge (1-\vartheta_j)\big(\psi_{j}^*- Q_\lambda(p_{j})\big) +\frac{1-\vartheta_j}{2\lambda\vartheta_j^2}\normt{p_{j+1} - q_j}^2 +  m_j(p_{j+1})\label{eq:analysis_TA_3}\\ %+ (1-\vartheta_j)(m_j(p_{j}) )
&\ge (1-\vartheta_j)\big(\psi_{j}^*- Q_\lambda(p_{j})\big)+ \delta_j +(L_j/2)\normt{p_{j+1} - q_j}^2  +  \phi(q_j) + \ipt{s_j}{p_{j+1}-q_j}\nn\\
&\hspace{4cm} + \omega_\lambda(p_{j+1}) - \delta_j, \quad \where\; \mbox{$\delta_j$ and $L_j$ are defined in~\eqref{eq:def_delta_L}}\label{eq:analysis_TA_4}\\
&\ge (1-\vartheta_j)\big(\psi_{j}^*- Q_\lambda(p_{j})\big)+ \phi(p_{j+1})+ \omega_\lambda(p_{j+1}) - \delta_j\label{eq:analysis_TA_5}\\
&= (1-\vartheta_j)\big(\psi_{j}^*- Q_\lambda(p_{j})\big)+ Q_\lambda(p_{j+1}) - \delta_j,
\end{align}
where~\eqref{eq:analysis_TA_1} follows from the $\lambda^{-1}$-strong convexity of $\psi_j$ on $\bbR^n$ and that $r_{j} = {\argmin}_{u\in\bbR^n}\;\psi_{j}(u)$ (cf.~\eqref{eq:def_psi_j}),~\eqref{eq:analysis_TA_2} follows from~\eqref{eq:triangle_identity},~\eqref{eq:analysis_TA_3} follows from the convexity of $m_j$ on $\bbR^n$,~\eqref{eq:analysis_TA_4} follows from Lemma~\ref{lem:lb_L_j}, and~\eqref{eq:analysis_TA_5} follows from~\eqref{eq:inexact_quad_ub}. 
As a result, we have 
\begin{align}
\Delta_{j+1}\le (1-\vartheta_j) \Delta_{j} +\delta_j = (1 - \vartheta_j/ 2)\Delta_{j} = \big(1 - (1/2)\min\{1/2, (\zeta_j/2)^{1/w}\}\big)\Delta_{j}. 
\end{align}
By the definition of $H$, we know that $j_0= \min\{0\le j \le J_\TA(\delta): (\zeta_j/2)^{1/w}\le  1/2\}$, and hence 
\begin{alignat}{2}
\Delta_{j+1}&\le (3/4) \Delta_{j}, &\qquad &\forall\, 0\le j< j_0, \label{eq:lin_rate_first}\\
\Delta_{j} - \Delta_{j+1}&\ge c\,\Delta_j^{\gamma+1}, &\qquad &\forall\,  j_0\le j <  J_\TA(\delta). \label{eq:sublin_second}
\end{alignat} 
 From~\eqref{eq:lin_rate_first}, we have 
\begin{equation}
\Delta_{j} \le (3/4)^j \Delta_0, \qquad\forall\, 0\le j\le j_0.
\end{equation}
Also, from~\eqref{eq:sublin_second}, by~\cite[Lemma 4.1]{Borwein_14}, we have 
\begin{align}
\Delta_j\le \big(\Delta_{j_0}^{-\gamma} + \gamma c(j-j_0)\big)^{-1/\gamma},\quad \forall\, j_0\le j\le  J_\TA(\delta).%\label{eq:rate_h_k}
\end{align} 
Since $\Delta_{j_0}\le H$, we complete the proof. 
\end{proof}

By Theorem~\ref{thm:TA}, we easily have the following corollary. For convenience, define
\begin{equation}
(x)_+:= \max\{x,0\}\quad \andd\quad  \ln_+(x) := \max\{\ln x,0\}. 
\end{equation}

\begin{corollary}\label{cor:complexity_TA}
In Algorithm~\ref{algo:triple-ave}, we have 
\begin{equation}
J_\TA(\delta)\le \left\lceil 4\ln_+\left(\frac{\Delta_0}{\max\{\delta,H\}}\right) \right\rceil + \left\lceil \frac{(\delta^{-\gamma} - \min\{\Delta_0,H\}^{-\gamma})_+}{\gamma c}\right\rceil,
\end{equation}
where  $H$, $\gamma$ and  $c$ are defined in Theorem~\ref{thm:TA}. As a result, we have 
\begin{align}
J_\TA(\delta) = O\left(\ln\left({D^2}{\lambda^{-1} H^{-1}}\right)+\delta^{-\gamma}\lambda^{\frac{1+\nu}{1+3\nu}}M^{\frac{2}{1+3\nu}}\right). 
\end{align} 
\end{corollary}

\subsection{Description and Analysis of Algorithm~\ref{algo:FW}} \label{sec:FW}

Algorithm~\ref{algo:FW} is the typical generalized FW method with ``open-loop'' stepsizes for solving smooth convex composite optimization problems (see e.g.,~\cite{Bach_15,Nest_18}). Fix any $j\ge 0$. %$0\le j\le J_{\FW}(\tau)$. 
Note that $G_j$ is precisely the (generalized) Wolfe-gap at $u_j$. In fact, it is well-known that $G_j$ is also the Fenchel duality gap at the primal-dual pair $(u_j,\nabla Q(u_j))$ (see e.g.,~\cite[Section~2]{Zhao_26a}), namely 
\begin{equation}
G_j = Q(u_j)+\barg(u_j) + Q^*(\nabla Q(u_j))+\barg^*(-\nabla Q(u_j)) \ge \barQ(u_j) - \barQ_*. \label{eq:duality_gap}
%\quad\andd\quad 
\end{equation}
We terminate Algorithm~\ref{algo:FW} as soon as $G_j\le \tau$. This criterion is justified by the following result. 

\begin{lemma} \label{lem:FW_output}
For any $j\ge 0$, we have
\begin{equation}
\barQ(u) - \barQ(u_j)\ge D_Q(u,u_j) - G_j\ge (2A)^{-1}\normt{u-u_j}^2 - G_j, \quad \forall\,u\in\bbR^n, \label{eq:quadratic_lb}
\end{equation}
where $D_Q$ denotes the Bregman divergence induced by $Q$, i.e., 
\begin{equation}
D_Q(u,u'): = \barQ(u) - \barQ(u') - \ip{\nabla \barQ(u')}{u-u'}, \quad \forall\,u,u'\in\bbR^n.  
\end{equation}
As a result, let $z$ be the output of Algorithm~\ref{algo:FW}, and we have 
\begin{equation}
\barQ(u) - \barQ(z)\ge (2A)^{-1}\normt{u-z}^2 - \tau, \quad \forall\,u\in\bbR^n. 
\end{equation}
\end{lemma}

\begin{proof}
By the optimality condition of the defining minimization problem of $v_j$ (cf.~Line~\ref{line:def_v_j}),
 we have 
 \begin{equation}
 \barg(u) +\ip{\nabla Q(u_j)}{u} \ge \barg(v_j) +\ip{\nabla Q(u_j)}{v_j}, \quad \forall\,u\in\bbR^n,
 \end{equation}
 which amounts to 
\begin{align}
&\barg(u) - \barg(u_j) +\ip{\nabla Q(u_j)}{u-u_j} \ge \barg(v_j) - \barg(u_j) +\ip{\nabla Q(u_j)}{v_j-u_j} = -G_j, \quad \forall\,u\in\bbR^n .
\end{align}
The first inequality in~\eqref{eq:quadratic_lb} follows from  
\begin{equation}
\begin{split}
\barg(u) - \barg(u_j) +\ip{\nabla Q(u_j)}{u-u_j} &= \barg(u) - \barg(u_j) -D_Q(u,u_j) - Q(u_j) + Q(u)\\
  &=  \barQ(u) - \barQ(u_j) - D_Q(u,u_j). 
\end{split}
\end{equation}
Also, the second inequality in~\eqref{eq:quadratic_lb} follows from the $A^{-1}$-strong convexity of $Q$ on $\bbR^n$.  %complete the proof, note that 
%\begin{equation}
%\end{equation}
\end{proof}

Next, let us analyze the total number of iterations of Algorithm~\ref{algo:FW}. 
%Let us also 
To that end, define the sequence of the ``best'' Wolfe-gaps $\{\uG_j\}_{j\ge 0}$ as 
\begin{equation}
\uG_j:= {\min}_{i=0,\ldots,j}\; G_i, \quad\forall\,j\ge 0. 
\end{equation}
Based on~\eqref{eq:duality_gap}, we can invoke~\cite[Proposition 3]{Bach_15} %we have the 
%which states following %result regarding the 
to obtain convergence rate of $\{\uG_j\}_{j\ge 0}$:

\begin{prop}\label{prop:bach}
In Algorithm~\ref{algo:FW}, %let $j$ be an index before 
 we have 
\begin{equation}
\uG_j \le \frac{8D^2}{A(j+2)}, \quad\forall\,j\ge 0. %\le j\le J_{\FW}(\tau). 
\end{equation}
\end{prop}

Let $J_{\FW}(\tau)$ be the iteration index at which Algorithm~\ref{algo:FW} terminates, % (at the accuracy level $\tau>0$), 
namely 
\begin{equation}
J_{\FW}(\tau):= \min\{j\ge 0: \, G_j\le \tau\}. 
\end{equation}
%Therefore, 
Based on Proposition~\ref{prop:bach}, we easily have the following corollary. 
\begin{corollary} \label{cor:comp_FW}
Define $\uJ_{\FW}(\tau):= \min\{j\ge 0: \, \uG_j\le \tau\}.$ Then
\begin{equation}
J_{\FW}(\tau) = \uJ_{\FW}(\tau) \le \left\lceil\frac{8D^2}{A\tau}\right\rceil -2. 
\end{equation}
\end{corollary}

\begin{proof}
By definition, it is easy to see that $J_{\FW}(\tau) = \uJ_{\FW}(\tau)$. Now, suppose that $\uJ_{\FW}(\tau) > j_0:= \left\lceil{8D^2}/({A\tau})\right\rceil -2$, then by Proposition~\ref{prop:bach}, we have $\uG_{j_0}\le \tau$, and hence $\uJ_{\FW}(\tau)\le j_0$, which leads to a contradiction. 
\end{proof}

\subsection{Analysis of Algorithm~\ref{algo:iAPG} } \label{sec:algo_iAPG}

%Fix any $\lambda>0$, and define $\Phi_\lambda^*:= {\min}_{u\in\bbR^n}\; \Phi_\lambda(u)$. 
%For any $x\in\bbR^n$, 

For any $k\ge 0$, define 
\begin{align}
\Psi_k(u)&:= \textstyle \sum_{i=0}^{k-1} \alpha_i\, \pi_i(u) + ({1}/{2}) \normt{u - a}^2,\\
&\hspace{1cm}\quad \where \quad \pi_i(u):= w_i + \ip{l_i}{u - y_i} + \barg(u),\quad \forall\,i\ge 0,\\ 
%\end{align}%Define%\begin{equation}
\barQ_k(u)&:= (2A_k)^{-1}\normt{u-c_k}^2 +\barg(u), \quad 
%\where\quad  %Q:= (2A)^{-1}\normt{\cdot-c}^2. 
%Q(u):= ,\;\;\; 
\forall\,u\in\bbR^n. 
\end{align}  
  %Proposition~\ref{prop:inexact_model_moreau_env}, 
Note that we have 
\begin{alignat}{3}
&\pi_i(u)\le\phi_\lambda(u) + \barg(u) = \Phi_\lambda(u),&\quad \forall\,u\in\bbR^n, &\quad \forall\,i\ge 0,\label{eq:pi_ub}\\
\quad \Longrightarrow\quad &\Psi_k(u)\le A_k \Phi_\lambda(u) +  ({1}/{2}) \normt{u - a}^2,&\quad \forall\,u\in\bbR^n, & \quad \forall\,k\ge 0.  \label{eq:Phi_ub}
\end{alignat}
Indeed, this follows from Lemma~\ref{lem:delta_L_inexact} when $\nu\in[0,1)$ and the convexity of $\phi$ when $\nu=1$. Also, note that %since $x_0 = a$, 
since $c_k = a-\xi_k$ for any $k\ge 1$, we can rewrite 
\begin{align}
\begin{split}
\Psi_k(u) = B_k   
&+ (1/2)\normt{u - c_k}^2 + A_k\barg(u)= B_k + A_k\barQ_k(u), \\ %(\barg(u) + (2A_k)^{-1}\normt{u - c_k}^2),\\
&\quad\where \quad B_k:= \textstyle \sum_{i=0}^{k-1} \alpha_i(w_i - \ip{l_i}{y_i})  - (1/2) \normt{\xi_k}^2 + \ip{\xi_k}{a}, \quad \forall\,k\ge 1.  
\end{split}\label{eq:Psi_k_rewrite}
\end{align}

%Define
%\begin{equation}
%K_{\rm APG}(\bar\delta) := \min\big\{k\ge 0: U_k\le \bar\delta\big\}. 
%\end{equation}
Now, let us present the convergence rate of Algorithm~\ref{algo:iAPG}. 

\begin{theorem} \label{thm:iAPG}
Fix any $\delta>0$. In Algorithm~\ref{algo:iAPG}, 
if $A_0 = 0$ and %for all $k\ge 0,$
\begin{equation}
A_{k+1}/\alpha_k^2\ge L, \quad\where\quad  L:=
\left\{\begin{array}{@{}ll}
2/\lambda, &\mbox{if }\; \nu\in[0,1)\\ % 
M, &\mbox{if }\; \nu=1
\end{array}
\right., \quad \forall\,k\ge 0,\label{eq:def_L}
\end{equation}
then we have 
\begin{align}
\gap_\lambda(x_{k}) = \Phi_\lambda(x_{k}) - \Phi_\lambda^*\le  \frac{\sum_{i=1}^{k} A_i\tau_i + 2\delta\sum_{i=1}^{k} A_i + D^2/2}{A_k}, \qquad \forall\,k\ge 1. %\le k\le K. 
\label{eq:Phi_lambda_gap} %_{\rm APG}(\bar\delta). 
\end{align}
\end{theorem}

\begin{proof}
From Algorithm~\ref{algo:triple-ave} and Lemma~\ref{lem:delta_L_inexact}, we know that if $\nu\in[0,1)$, then %for all $k\ge 0$, 
\begin{align}
0\le \phi_\lambda(x)  - w_k - \ip{l_k}{x-y_k} \le \lambda^{-1} \normt{x-y_k}^2 + 2\delta,  \quad \forall\, x\in\bbR^n, \quad \forall\, k\ge 0. \label{eq:Moreau_consequence}
\end{align}
If $\nu=1$, then $(w_k,l_k) = (\phi(y_k),\nabla \phi(y_k))$, and hence 
\begin{align}
0\le \phi_\lambda(x)  - w_k - \ip{l_k}{x-y_k} \le (M/2) \normt{x-y_k}^2,  \quad \forall\, x\in\bbR^n, \quad \forall\, k\ge 0.  \label{eq:smoothness_phi}
\end{align}
Therefore, by the definition of $L$ in~\eqref{eq:def_L},  for all  $\nu\in[0,1]$, we have
\begin{align}
0\le \phi_\lambda(x)  - w_k - \ip{l_k}{x-y_k} \le (L/2) \normt{x-y_k}^2+ 2\delta,  \quad \forall\, x\in\bbR^n, \quad \forall\, k\ge 0.  \label{eq:inexact_moreau_combined}
\end{align}
From Lemma~\ref{lem:FW_output} and~\eqref{eq:Psi_k_rewrite},  we know that for all $k\ge 1$, 
\begin{align}
&\barQ_k(u) \ge  \barQ_k(z_k)+ (2A_k)^{-1}\normt{u-z_k}^2 - \tau_k, \quad \forall\,u\in\bbR^n.\label{eq:FW_conseq_0}
\end{align}
%Consequently, this implies that 
%By~\eqref{eq:Psi_k_rewrite}, %and the fact that $\Psi_0(u) = (1/2)\normt{u-a}^2$, $A_0=0$ and $z_0 = a$, %hence $\Psi_0(z_0) = (1/2)\normt{z_0-a}^2 = (1/2)\normt{a-a}^2=0$, 
%we have %implies that %for all $k\ge 0$, 
Next, note that we have
\begin{align}
&\Psi_k(u) \ge  \Psi_k(z_k) + (1/2)\normt{u-z_k}^2 - A_k\tau_k, \quad \forall\,u\in\bbR^n,\quad \forall\,k\ge 0. \label{eq:FW_consequence}
\end{align}
%Note that since 
Indeed, for $k=0$, this follows from $\Psi_0(u) = (1/2)\normt{u-a}^2$, $A_0=0$ and $z_0 = a$, %$A_0 = 0$ and $z_0 = a$, 
which then implies that $\Psi_0(z_0)=0$; for $k\ge 1$, this follows from~\eqref{eq:FW_conseq_0} and~\eqref{eq:Psi_k_rewrite}.  %  As such,~\eqref{eq:FW_consequence} holds for $k=0$ as well. 
Let $\barx^*$ be an optimal solution of~\eqref{eq:prob_Phi_lambda}. 
%This, 
By taking $u=\barx^*$ in~\eqref{eq:FW_consequence}, and using~\eqref{eq:Phi_ub} and the fact that $\barx^*,a\in %\dom\barg = 
\dom g$, 
we have %implies that  
\begin{align}
%\Longrightarrow\qquad & 
A_k \Phi_\lambda^* +   D^2/2 \ge A_k \Phi_\lambda(\barx^*) +  ({1}/{2}) \normt{\barx^* - a}^2 \ge \Psi_k(\barx^*) \ge  \Psi_k(z_k) %+ (1/2)\normt{\barx^*-z_k}^2 
- A_k\tau_k, \quad \forall\,k\ge 0. \label{eq:ub_Psi_k_z_k}  
\end{align}
Now, for any $k\ge 0$, we have %that for any $u\in\bbR^n$, 
\begin{align}
\Psi_{k+1}(z_{k+1}) &= \Psi_{k}(z_{k+1}) + \alpha_k \pi_k(z_{k+1})\\
& \ge  \Psi_k(z_k) + (1/2)\normt{z_{k+1}-z_k}^2 - A_k\tau_k + \alpha_k \pi_k(z_{k+1}) \label{eq:iAPG_1}\\
&\ge  \Psi_k(z_k) + (2\theta_k^2)^{-1}\normt{x_{k+1}-y_k}^2 - A_k\tau_k + \alpha_k \pi_k(z_{k+1}) + A_k (\pi_k(x_k) - \Phi_\lambda(x_k))\label{eq:iAPG_2}\\
&\ge \Psi_k(z_k) + (2\theta_k^2)^{-1}\normt{x_{k+1}-y_k}^2 + A_{k+1} \pi_k(x_{k+1})- A_k \Phi_\lambda(x_k) - A_k\tau_k\label{eq:iAPG_3}\\
%&\hspace{5cm} %+ A_k \pi_k(x_k) 
&= \Psi_k(z_k) + A_{k+1}\big\{2\delta+(2\theta_k^2A_{k+1})^{-1}\normt{x_{k+1}-y_k}^2 +  w_k + \ip{l_k}{x_{k+1} - y_k} + \barg(x_{k+1})\big\}\nn\\
&\hspace{6cm} - A_k \Phi_\lambda(x_k) - A_k\tau_k - 2\delta A_{k+1} \\ %\label{eq:iAPG_4}\\
&\ge \Psi_k(z_k) + A_{k+1}\Phi_\lambda(x_{k+1}) %\big\{\phi_\lambda(x_{k+1}) + \barg(x_{k+1})\big\}
 - A_k \Phi_\lambda(x_k) - A_k\tau_k - 2\delta A_{k+1}, \label{eq:iAPG_4}
\end{align}
where in~\eqref{eq:iAPG_1} we use~\eqref{eq:FW_consequence}, in~\eqref{eq:iAPG_2} we use $z_{k+1} - z_k = \theta_k(x_{k+1} - y_k)$ and~\eqref{eq:pi_ub}, in~\eqref{eq:iAPG_3} we use %the definition of 
$\theta_k = \alpha_k/A_{k+1}$ (cf.~\eqref{eq:def_alpha_theta}) and the convexity of $\pi_k$, in~\eqref{eq:iAPG_4} we use~\eqref{eq:def_L} and~\eqref{eq:inexact_moreau_combined}.

Next, fix any $k\ge 1$ and telescope~\eqref{eq:iAPG_4} over $i = 0,\ldots\,k-1$, we have
\begin{align}
\Psi_{k}(z_{k})\ge \Psi_0(z_0) + A_{k}\Phi_\lambda(x_{k})  - A_0 \Phi_\lambda(x_0) - \textstyle\sum_{i=0}^{k-1} A_i\tau_i - 2\delta \textstyle\sum_{i=1}^{k} A_{i},\label{eq:telescoped_algo_iAPG}
\end{align}
By rearranging~\eqref{eq:telescoped_algo_iAPG}, and noting that $A_0=0$ and $z_0 = a$ and $\Psi_0(z_0)=0$, we have 
\begin{align}
\textstyle\sum_{i=1}^{k-1} A_i\tau_i + 2\delta\sum_{i=1}^{k} A_i &=\textstyle\sum_{i=0}^{k-1} A_i\tau_i + 2\delta\sum_{i=1}^{k} A_i + A_{0}\Phi_\lambda(x_{0}) - \Psi_0(z_0)\\
& \ge  A_{k}\Phi_\lambda(x_{k}) - \Psi_{k}(z_{k}) \\
&\ge  A_{k}\Phi_\lambda(x_{k}) - A_k \Phi_\lambda^* -   D^2/2  - A_k\tau_k,
\end{align}
where the last inequality follows from~\eqref{eq:ub_Psi_k_z_k}.  
%Since $A_0 = 0$, 
We then reach~\eqref{eq:Phi_lambda_gap}. %Therefore,
%\begin{align}
%\Phi_\lambda(x_{k}) - \Phi_\lambda^*\le \frac{\sum_{i=1}^{k} A_i\tau_i + 2\delta\sum_{i=1}^{k} A_i + D^2/2}{A_k}. \tag*{$\square$} %  = U_k. 
%\end{align}
%\renewcommand{\qed}{}
\end{proof}

Based on Theorem~\ref{thm:iAPG}, we easily have the following corollaries. 

\begin{corollary}\label{cor:bound_gap_lambda_x_K}
%In particular, 
Fix any $\lambda,\bar\delta>0$ and let $L$ be defined in~\eqref{eq:def_L}. % Theorem~\ref{thm:iAPG}. 
In Algorithm~\ref{algo:iAPG}, if we choose % $A_k:= k^2/L$, $ $then 
%$\lambda := 1/(18L_\varepsilon),$ %$\bar\delta:=(2/9)\varepsilon, $
\begin{equation}
\begin{split}
&K := \left\lceil\sqrt{8 L D^2/\bar\delta}\right\rceil,\quad A_k := k^2/(4L),\quad \forall\,0\le k\le K,  \\
&%\qquad\qquad
\delta := \frac{\bar\delta}{2(K+1)}\quad\andd\quad \tau_k = \frac{\bar\delta K}{16LA_k}, \quad \forall\,1\le k\le K, 
\end{split}\label{eq:iAPG_param_choices}
\end{equation} 
then $\gap_\lambda(x_{K}) %= \Phi_\lambda(x_{K}) - \Phi_\lambda^*
\le \bar\delta$. 
%\begin{align}
%\delta
%\end{align}
\end{corollary}

\begin{proof}
First, note that by the choice of $\{A_k\}_{0\le k\le K}$ in~\eqref{eq:iAPG_param_choices}, we indeed have $A_0 = 0$ and 
\begin{equation}
A_{k+1}/\alpha_k^2 = A_{k+1}/(A_{k+1} - A_k)^2 = 4L (k+1)^2/(2k+1)^2\ge L, \quad\forall\,0\le k\le K-1. 
\end{equation}
Also, note that by the choices of the parameters in~\eqref{eq:iAPG_param_choices}, we have  
\begin{align}
&\frac{D^2}{2A_K} = \frac{2LD^2}{K^2} \le  \frac{\bar\delta}{4},\\
&2\delta\frac{\sum_{i=1}^{K} A_i}{A_K} = 2\delta\frac{\sum_{i=1}^{K} i^2}{K^2} = \delta\frac{(K+1)(2K+1)}{3K}\le (K+1)\delta= \frac{\bar\delta}{2},\\
&\frac{\sum_{i=1}^{K} A_i\tau_i}{A_K} = \frac{\bar\delta K^2/(16L)}{K^2/(4L)} = \frac{\bar\delta}{4}. 
\end{align}
These bounds, together with~\eqref{eq:Phi_lambda_gap}, complete the proof.  
\end{proof}

%Based on Theorem~\ref{thm:iAPG}, we easily have the following corollary. 

\begin{corollary}
Fix any $\varepsilon>0$.   
In Algorithm~\ref{algo:iAPG},  if we choose $\lambda := 1/(18L_\varepsilon)$, 
 $\bar\delta := (2/9)\varepsilon$, and the other parameters according to~\eqref{eq:iAPG_param_choices},  then $%P(x_{K}) - P^*
\gap(x_{K})\le \varepsilon/2$. 
\end{corollary}

\begin{proof}
By Corollary~\ref{cor:bound_gap_lambda_x_K}, we know that $\gap_\lambda(x_{K})\le\bar\delta= (2/9)\varepsilon$. 
We then invoke Lemma~\ref{lem:bound_gap} to complete the proof. Specifically, the choice of $\lambda$ %in~\eqref{eq:iAPG_param_choices} 
corresponds to choosing $c = 1/9$ in~\eqref{eq:bound_gap}, and by letting $y = x_K$ in~\eqref{eq:bound_gap}, we have % which then becomes %and hence
\begin{equation}
\gap(x_K)\le \frac{(2/9) \varepsilon  + \gap_\lambda(x_K)}{8/9}\le \frac{(2/9) \varepsilon  + (2/9)\varepsilon}{8/9} = \frac{\varepsilon}{2}. \tag*{$\square$} %\label{eq:bound_gap}
\end{equation}
\renewcommand{\qed}{}
\end{proof}

\begin{theorem} \label{thm:comp_iAPG}
Fix any $\varepsilon>0$.
In Algorithm~\ref{algo:iAPG},  choose $\lambda := 1/(18L_\varepsilon)$, 
 $\bar\delta := (2/9)\varepsilon$, and the other parameters according to~\eqref{eq:iAPG_param_choices}. Let $N_{\rm GLMO}$ and $N_{\rm FO}$ denote the number of GLMO calls and FO calls in Algorithm~\ref{algo:iAPG}, respectively. Then for all $\nu\in[0,1]$, we have 
\begin{align}
N_{\rm GLMO} = O\left(   M^{\frac{2}{1+\nu}}D^2 \varepsilon^{-\frac{2}{1+\nu}}\right) \quad \andd \quad N_{\rm FO} = O\Big(M^{\frac{2}{1+3\nu}}D^{{\frac{2(1+\nu)}{1+3\nu}}}\varepsilon^{-\frac{2}{1+3\nu}}\Big).
\end{align}
\end{theorem}

\begin{proof}
By the definitions of $L$ (cf.~\eqref{eq:def_L}) and $\lambda$, we have $L =\Theta(\lambda^{-1})= \Theta(L_\varepsilon)$ for $\nu\in[0,1)$ and $L= \Theta(M) = \Theta(L_\varepsilon)$ for $\nu=1$. Thus $L= \Theta(L_\varepsilon)$ for all $\nu\in[0,1]$. Also, we clearly have $\bar\delta = \Theta(\varepsilon)$. In addition, each iteration of Algorithm~\ref{algo:FW} makes one GLMO call. As a result, from Corollary~\ref{cor:comp_FW}, we know that for all $\nu\in[0,1]$, 
\begin{align}
&N_{\rm GLMO}= O\left(  \sum_{i=1}^{K} \frac{D^2}{A_i\tau_i}\right) = O\left(  \sum_{i=1}^{K} \frac{LD^2}{\bar\delta K}\right) = O\left(   \frac{L_\varepsilon D^2}{\varepsilon}\right) 
= O\left(   M^{\frac{2}{1+\nu}}D^2 \varepsilon^{-\frac{2}{1+\nu}}\right).
\end{align}
Next, let $H$, $\gamma$ and  $c$ be defined in Theorem~\ref{thm:TA}. Note that each iteration of Algorithm~\ref{algo:triple-ave} makes one FO call. From Corollary~\ref{cor:complexity_TA} and the definition of $c$, we know that for $\nu\in[0,1)$, 
\begin{align}
&N_{\rm FO} = O\left(  K\big(\ln\left({D^2}{\lambda^{-1} H^{-1}}\right)+\delta^{-\gamma} c^{-1}
%\lambda^{\frac{1+\nu}{1+3\nu}}M^{\frac{2}{1+3\nu}}
\big)\right).
\end{align}
Note that since $H = \Theta(\lambda^{\frac{1+\nu}{1-\nu}}M^{\frac{2}{1-\nu}}) = \Theta(L_\varepsilon^{-\frac{1+\nu}{1-\nu}}M^{\frac{2}{1-\nu}})=\Theta(\varepsilon)$, %and $\bar\delta =\Theta(\varepsilon)$, 
we have 
\begin{align}
{D^2}{\lambda^{-1} H^{-1}} = \Theta(LD^2/\varepsilon) = \Theta(K^2)\quad\Longrightarrow\quad \ln\left({D^2}{\lambda^{-1} H^{-1}}\right) = \Theta(\ln K).
\end{align}
Also, note that $\delta = \Theta(\varepsilon/K)$ and $c^{-1} = \Theta(\lambda^{\frac{1+\nu}{1+3\nu}}M^{\frac{2}{1+3\nu}}) = \Theta(L_\varepsilon^{-\frac{1+\nu}{1+3\nu}}M^{\frac{2}{1+3\nu}}) = \Theta(\varepsilon^{\frac{1-\nu}{1+3\nu}}) = \Theta(\varepsilon^\gamma)$, and hence $\delta^{-\gamma} c^{-1}=\Theta(K^\gamma)$. %As a result, 
Since $\nu\in[0,1)$, we have $\gamma>0$ and hence 
\begin{align}
N_{\rm FO} = O\left(  K(\ln K+ K^\gamma)\right) &= O\left( K^{\gamma+1}\right) = O\Big(\big(L_\varepsilon D^2/\varepsilon\big)^{\frac{\gamma+1}{2}}\Big)\nn\\
& = O\left(\big(\varepsilon^{-\frac{2}{1+\nu}}D^2M^{\frac{2}{1+\nu}}\big)^{\frac{1+\nu}{1+3\nu}}\right) = O\Big(M^{\frac{2}{1+3\nu}}D^{{\frac{2(1+\nu)}{1+3\nu}}}\varepsilon^{-\frac{2}{1+3\nu}}\Big). \label{eq:comp_nu<1} %\right). 
\end{align}
Now, if $\nu=1$, then %$L = M$ and 
each call of Algorithm~\ref{algo:triple-ave} within Algorithm~\ref{algo:iAPG} only makes one FO call, and hence 
\begin{align}
N_{\rm FO} = O(K) = O\big(M^{1/2} D \varepsilon^{-1/2}\big). 
\end{align} 
Note that this complexity precisely corresponds to~\eqref{eq:comp_nu<1} for $\nu =1$. 
\end{proof}

\subsection{Total LO and GLMO Complexities of Algorithm~\ref{algo:restart}} \label{sec:total_comp}

Based on Theorem~\ref{thm:comp_iAPG}, we can easily analyze the total LO and GLMO complexities of Algorithm~\ref{algo:restart}. 

\begin{theorem}
%needed to produce $$ In addition,
Fix any $0<\bar\varepsilon< \varepsilon_0$. In Algorithm~\ref{algo:restart},  let $S= S(\bar\varepsilon): =\left\lceil \log_2(\varepsilon_0/\bar\varepsilon)\right\rceil$, then we have $\gap(a_S)\le\bar\varepsilon$. %we have 
%\begin{equation}
%\gap(a_S)\le \varepsilon_s = 2^{-s} \varepsilon_0. 
%\end{equation}
%In addition, let 
In addition, let $N_{\rm GLMO}(S)$ and $N_{\rm FO}(S)$ denotes the number of GLMO and FO calls  made by Algorithm~\ref{algo:restart} to produce $a_{S}$, then we have 
\begin{equation}
N_{\rm GLMO}(S)= O\left(   M^{\frac{2}{1+\nu}}D^2 \bar\varepsilon^{-\frac{2}{1+\nu}}\right) \quad \andd \quad N_{\rm FO}(S) = O\Big(M^{\frac{2}{1+3\nu}}D^{{\frac{2(1+\nu)}{1+3\nu}}}\bar\varepsilon^{-\frac{2}{1+3\nu}}\Big).
\end{equation}
\end{theorem}

\begin{proof}
Note that in Algorithm~\ref{algo:restart}, for any $s\ge 0$, we have 
\begin{equation}
\gap(a_s)\le \varepsilon_s = 2^{-s} \varepsilon_0. 
\end{equation}
%Without loss of generality, assume $\bar\varepsilon< \varepsilon_0$.
% Define $S:= \min\{s\ge 0: \gap(a_s)\le \bar\varepsilon\}$. Clearly, we have 
% $S \le  \left\lceil \log_2(\varepsilon_0/\bar\varepsilon)\right\rceil$. % = O(\ln(\varepsilon_0/\bar\varepsilon))$
%\begin{equation}
%By the definition of 
Since $S= S(\bar\varepsilon)$, we clearly have $\gap(a_S)\le\bar\varepsilon$. 
%S \le  \left\lceil 2\ln_+(\varepsilon_0/\bar\varepsilon)\right\rceil
%\end{equation}
Next, note that $N_{\rm GLMO}(S)$ and $N_{\rm FO}(S)$ are equal to the total number of GLMO calls and FO calls  made by Algorithm~\ref{algo:restart}  after running iteration $S-1$. % iterations of Algorithm~\ref{algo:restart}, 
Then from Theorem~\ref{thm:comp_iAPG}, we have that for all $\nu\in[0,1]$, 
\begin{align}
&N_{\rm GLMO} = O\left( \textstyle \sum_{s=0}^{S-1}  M^{\frac{2}{1+\nu}}D^2 \varepsilon_s^{-\frac{2}{1+\nu}}\right) = O\left( \textstyle   M^{\frac{2}{1+\nu}}D^2 \varepsilon_0^{-\frac{2}{1+\nu}} \sum_{s=0}^{S-1}2^{\frac{2s}{1+\nu}}\right)\nn \\
&\quad= O\left( \textstyle   M^{\frac{2}{1+\nu}}D^2 \varepsilon_0^{-\frac{2}{1+\nu}} 2^{\frac{2S}{1+\nu}}\right) = O\left( \textstyle   M^{\frac{2}{1+\nu}}D^2 \varepsilon_0^{-\frac{2}{1+\nu}} (\varepsilon_0/\bar\varepsilon)^{\frac{2}{1+\nu}}\right)  = O\left( \textstyle   M^{\frac{2}{1+\nu}}D^2  \bar\varepsilon^{-\frac{2}{1+\nu}}\right), 
\end{align}
and using the same reasoning, we have
 \begin{align}
N_{\rm FO} = O\Big(\textstyle \sum_{s=0}^{S-1}  M^{\frac{2}{1+3\nu}}D^{{\frac{2(1+\nu)}{1+3\nu}}}\varepsilon_s^{-\frac{2}{1+3\nu}}\Big) = O\Big(\textstyle  M^{\frac{2}{1+3\nu}}D^{{\frac{2(1+\nu)}{1+3\nu}}}\bar\varepsilon^{-\frac{2}{1+3\nu}}\Big). \tag*{$\square$}
\end{align}
\renewcommand{\qed}{}
\end{proof}

\paragraph{AI Use Disclosure.}
The author used ChatGPT 6 (Sol, Ultra mode) %. The AI tools were used inainly for  
for  literature search and exploring and stress-testing proof strategies, during the preparation of this manuscript. The author takes full responsibility for the mathematical claims, citations, and final text.

%\newpage
\bibliographystyle{abbrv}
\bibliography{math_opt}

\end{document}